\documentclass[12pt]{amsart}

\usepackage[T1]{fontenc}
\usepackage{lmodern}
\usepackage{amsmath,amssymb,amsfonts,mathrsfs,mathtools}
\usepackage{array,longtable,booktabs}
\usepackage{enumitem}
\usepackage{microtype}
\usepackage{graphicx}
\usepackage{xcolor}
\usepackage{needspace}
\usepackage{tikz-cd}
\usepackage[
  colorlinks=true,
  linkcolor=red,
  citecolor=red,
  urlcolor=red,
  filecolor=red,
  pdfpagemode=UseOutlines,
  bookmarks=true,
  bookmarksopen=true,
  pagebackref=true
]{hyperref}
\usepackage{bookmark}

\hypersetup{
  pdftitle={Concatenation perturbations and classical cumulants},
  pdfauthor={},
  pdfsubject={Double tensor unshuffle Hopf algebras, fixed-point perturbations, Hopf quotients, and comodule structure},
  pdfkeywords={double tensor algebra, unshuffle bialgebra, half-shuffle fixed-point equation, concatenation perturbation, free cumulants, classical cumulants, non-crossing closure, Hopf quotient},
  pdfstartview=FitH
}

\allowdisplaybreaks
\setlist[enumerate]{leftmargin=*,label=\textup{(\arabic*)}}

\newtheorem{theorem}{Theorem}[section]
\newtheorem{proposition}[theorem]{Proposition}
\newtheorem{lemma}[theorem]{Lemma}
\newtheorem{corollary}[theorem]{Corollary}
\theoremstyle{definition}
\newtheorem{definition}[theorem]{Definition}
\newtheorem{example}[theorem]{Example}
\newtheorem{remark}[theorem]{Remark}

\def\shascale{0.75}
\makeatletter
\newcommand{\qshu@scaledsha}[2]{%
  \vcenter{\hbox{\scalebox{\shascale}{$#1#2$}}}%
}
\IfFileExists{shuffle.sty}{%
  \usepackage{shuffle}
  \DeclareRobustCommand{\sha}{%
    \ensuremath{\mathbin{\mathpalette\qshu@scaledsha\shuffle}}%
  }%
}{%
  \DeclareRobustCommand{\sha}{%
    \ensuremath{\mathbin{\mathpalette\qshu@scaledsha\amalg}}%
  }%
}
\makeatother

\newcommand{\A}{\mathcal A}
\newcommand{\deshat}{\Delta_{\sha}^{\overline{T}}}
\newcommand{\deshah}{\Delta_{\sha}^{H(\A)}}
\newcommand{\reddeshat}{\bar{\Delta}_{\sha}^{\overline{T}}}
\newcommand{\reddeshah}{\bar{\Delta}_{\sha}^{H(\A)}}
\newcommand{\clNC}{\operatorname{cl}_{\NC}}
\newcommand{\Ext}{\operatorname{Ex}}
\newcommand{\Lin}{\operatorname{Lin}}

\newcommand{\onepart}[1]{\mathbf 1_{#1}}
\newcommand{\Conn}{\mathrm{Conn}}
\newcommand{\conc}{\mathfrak{c}}
\newcommand{\sesec}{\mathfrak{b}}
\newcommand{\NC}{\mathrm{NC}}
\newcommand{\Part}{\mathcal P}
\newcommand{\ioconc}{{\iota^{\mathfrak{c}}}}
\newcommand{\Phiconc}{\Phi^{\mathfrak{c}}}
\newcommand{\cconc}{c_H^{\mathfrak{c}}}
\newcommand{\rhoR}{\rho_{\mathrm R}}
\newcommand{\rhoL}{\rho_{\mathrm L}}
\newcommand{\idH}{\operatorname{id}_{H(\A)}}
\newcommand{\idTbar}{\operatorname{id}_{\overline T(\A)}}
\newcommand{\DprecH}{\Delta^{+}_{\prec,H(\A)}}
\newcommand{\DsuccH}{\Delta^{+}_{\succ,H(\A)}}
\newcommand{\redDprecH}{\Delta_{\prec,H(\A)}}
\newcommand{\redDsuccH}{\Delta_{\succ,H(\A)}}
\newcommand{\DprecTbar}{\Delta^{+}_{\prec,\overline T(\A)}}
\newcommand{\DsuccTbar}{\Delta^{+}_{\succ,\overline T(\A)}}
\newcommand{\redDprecTbar}{\Delta_{\prec,\overline T(\A)}}
\newcommand{\redDsuccTbar}{\Delta_{\succ,\overline T(\A)}}
\newcommand{\convTbar}{*_{\overline T(\A)}}
\newcommand{\convH}{*_{H(\A)}}
\newcommand{\bfk}{{\bf k}}
\newcommand{\mybullet}{\mathbin{\vcenter{\hbox{\scalebox{0.5}{$\bullet$}}}}}
\providecommand{\Item}{{\rm Item}}

\newcommand{\nc}{\newcommand}
\nc{\li}[1]{\textcolor{red}{#1}}
\nc{\lir}[1]{\textcolor{red}{Li:#1}}
\nc{\yi}[1]{\textcolor{blue}{Yi: #1}}
\nc{\xing}[1]{\textcolor{purple}{Xing:#1}}

\begin{document}

\title[Classical cumulants from the double tensor unshuffle Hopf algebra]{Classical cumulants from the double tensor unshuffle Hopf algebra}

\author{Xing Gao}
\address{School of Mathematics and Statistics, Lanzhou University,
Lanzhou, 730000, China; Gansu Provincial Research Center for Basic Disciplines of Mathematics and Statistics, Lanzhou, 730070, China
}
\email{gaoxing@lzu.edu.cn}

\author{Jingtao Li}
\address{School of Mathematics and Statistics, Lanzhou University
Lanzhou, 730000, China
}
\email{ljingtao2025@lzu.edu.cn}

\date{\today}

\begin{abstract}
We realize classical cumulants and their fixed-point equations in the
double tensor unshuffle Hopf algebra introduced by Ebrahimi-Fard and
Patras, placing classical, free, Boolean, and monotone cumulants in a
common algebraic framework. A concatenation perturbation of the
functionals yields left and right half-shuffle equations on ordinary
words without changing the coproduct or either half-coproduct.
Both equations are equivalent to the classical moment--cumulant
formula and uniquely determine the classical cumulants.
We also obtain a Hopf-algebraic realization of the classical-to-free
cumulant transform through convolution exponentiation and the left
half-shuffle logarithm. A commutative diagram identifies its
partition-level action with summation over non-crossing closure
fibers, recovering Lehner's closure-fiber formula within the
framework of the double tensor unshuffle Hopf algebra.
\end{abstract}


\maketitle
\begingroup
\small\linespread{0.95}\selectfont
\tableofcontents
\endgroup

\section{Introduction}
\label{sec:introduction}

\subsection{Background and motivation}
\label{subsec:intro-motivation}

Cumulants describe moment functionals in terms adapted to different
notions of independence. The corresponding moment--cumulant formulas
are governed by different partition structures: all set partitions
in the classical case, non-crossing partitions in the free case,
interval partitions in the Boolean case, and non-crossing partitions
with additional block-ordering data in the monotone
case~\cite{AHLV2015,HS2011,NS2006,Speed1983,Speicher1994,SW1997}. A common algebraic description
of these formulas can therefore explain both their similarities and
the transformations relating their cumulants.

The double tensor unshuffle Hopf algebra $H(\A)$ developed by
Ebrahimi-Fard and Patras provides such a description for free, Boolean,
and monotone cumulants~\cite{CEFP2022,EFP2015,EFP2018,EFP2019,EFP2020}. In this setting,
moments are encoded by characters and cumulants by infinitesimal
characters. The left and right half-shuffle exponentials describe
the free and Boolean moment--cumulant relations, respectively,
whereas the convolution exponential describes the monotone relation.
Their inverse logarithms recover cumulants from moments, making
these three theories accessible through operations on the same
Hopf algebra.

Classical cumulants also admit a half-shuffle description, but on
the ordinary tensor algebra equipped with its cocommutative
unshuffle coproduct~\cite[Section~7]{EFP2015}. This leads to a
natural question: can classical cumulants and their fixed-point
equations be realized on $H(\A)$, the same Hopf algebra used for
the other three theories? Our aim is to retain
its coproduct and half-coproducts, so that all four
moment--cumulant correspondences can be studied within a single
ambient algebraic framework.

The relevant distinction is visible in the recursion at the block
containing the first index. The free recursion multiplies the
moments of the separate complement gaps, whereas the classical
recursion takes the moment of the entire ordered
complement~\cite[proof of Theorem~5.2 and Section~7]{EFP2015}.
The coproduct of $H(\A)$ records these gaps as separate bar
components. Deleting the bars merges the components while preserving
the order of the letters. This suggests perturbing the functionals
that evaluate the half-coproducts rather than modifying the
underlying unshuffle structure.

A second motivation concerns the passage from classical to free
cumulants. Non-crossing closure sends a set partition to its least
non-crossing coarsening. Lehner's formula describes the corresponding
cumulant transform by summation over the fibers of this map; in
particular, free cumulants are sums of classical cumulant products
over connected partitions~\cite[Theorem~4.1 and its proof]{Lehner2002}.
The multivariate relation is recorded in
\cite[Section~1]{AHLV2015}. We seek to describe this combinatorial transform in terms of the
Hopf-algebraic operations on $H(\A)$ and connect it with the
perturbed half-shuffle equations.

\subsection{Fixed-point equations and cumulant transformations}
\label{subsec:intro-contributions}

Our results extend both the fixed-point description of cumulants
and the transformations between them while keeping the double tensor
unshuffle Hopf algebra $H(\A)$ fixed. In the approach of
Ebrahimi-Fard and Patras, the free and Boolean cumulant infinitesimal
characters $\kappa$ and $\beta$ associated with a moment character
$\Phi$ are characterized by
\[
\Phi=\varepsilon_{H(\A)}+\kappa\prec\Phi,
\qquad
\Phi=\varepsilon_{H(\A)}+\Phi\succ\beta,
\]
respectively. These identities hold on all of $H(\A)$, and their
solutions are given by the left and right half-shuffle
exponentials~\cite{EFP2015,EFP2018,EFP2019}.

For the monotone cumulant infinitesimal character $\rho$, the
standard description is instead the convolution-exponential identity
$\Phi=\exp^{\convH}(\rho)$. The pre-Lie Magnus expansion and its
compositional inverse $W$ relate this identity to half-shuffle
fixed-point equations through
$\kappa=W(\rho)$ and $\beta=-W(-\rho)$~\cite{EFP2018}.
The coefficients in these equations are therefore transformed
cumulants rather than $\rho$ itself. This distinction shows that
a common Hopf-algebraic description need not take the form of the
same unperturbed half-shuffle equation for every cumulant family.

The present work brings classical cumulants into this picture
through the perturbed fixed-point equations
\eqref{eq:classical-left-fixed-point} and~\eqref{eq:right-fixed-point} below.
They modify the moment functional in the moment-perturbed fixed-point equation
and the cumulant functional in the right one, respectively.
Theorem~\ref{thm:classical-perturbations} shows that both equations
are equivalent to the classical moment--cumulant formula and uniquely
determine the classical cumulants from the moments.
Their significance lies in making the existing half-shuffle
operations applicable to classical cumulants without replacing
the coproduct or either half-coproduct. The distinction between
the classical and noncommutative cumulant descriptions is thus
expressed at the level of functional evaluation within a fixed
unshuffle structure. Unlike the standard free and Boolean
identities, the classical equations are imposed only on ordinary
words; moreover, the perturbed moment and cumulant functionals
need not remain characters and infinitesimal characters,
respectively.

The Hopf-quotient interpretation of bar deletion gives this
comparison structural content beyond a reformulation of scalar
recursions. Theorem~\ref{thm:flattening-morphism} shows that the
bar-deleting map realizes $\overline T(\A)$ as a Hopf quotient of
$H(\A)$ and preserves both half-coproducts, thereby connecting
the tensor-algebra description of classical cumulants with the
operations used for the other three families.
Combinatorially, the perturbation replaces the product of moments
over separate complement gaps by the moment of their concatenation.
In the first-block recursion, the resulting correction consists
precisely of the cumulant products indexed by complement partitions
with a block crossing the distinguished block, as expressed in
\eqref{eq:local-crossing-correction} and
Lemma~\ref{lem:gap-crossing}. Thus the same unshuffle structure
accommodates different partition combinatorics through different
functional evaluations, without identifying the underlying
notions of independence.

This perspective also connects the perturbed classical equations
with the standard free equation through their common moment
character. Starting from a classical cumulant word functional $c$,
our construction recovers the moment character $\Phi_c$ and then
applies the left half-shuffle logarithm to obtain
$\mathcal L(c)=\log_{\prec}^{H(\A)}(\Phi_c)$.
The commutative diagram of Theorem~\ref{thm:expanded-diagram}
identifies this algebraic passage with summation over non-crossing
closure fibers, recovering Lehner's known
formula~\cite{AHLV2015,Lehner2002} within this Hopf-algebraic
construction. It therefore relates two cumulant descriptions of
the same moment data and connects the fixed-point characterization
with the connected-partition expansion.

The significance of this construction extends to the mutual
conversion of all four cumulant families. Free, Boolean, and
monotone cumulants already admit invertible transformations
obtained from half-shuffle and convolution exponentials and
logarithms, equivalently expressed using the pre-Lie Magnus
expansion and its inverse~\cite{EFP2018,EFP2019}.
Our realization of the classical-to-free transform connects
classical cumulants to these existing transformations.
By Lemma~\ref{lemma:closure-fiber-formula},
$\mathcal L(c)(a_I)$ equals $c(a_I)$ plus products of cumulants
of strictly shorter subwords, so $\mathcal L$ is invertible
recursively in word length. Composing $\mathcal L$ and its inverse
with the established transformations therefore gives conversions
in both directions between any two of the classical, free,
Boolean, and monotone cumulant families. In this sense, the present
construction completes their mutual conversion within a common
Hopf-algebraic formalism, with all four families referring to the
same moment character.

\subsection{Outline of the paper}
Section~\ref{sec:perturbed-equations} introduces the concatenation
perturbation and the associated functionals
(Definition~\ref{def:perturbed-functionals}). Its main result is
Theorem~\ref{thm:classical-perturbations}, which characterizes
classical cumulants by either of the two perturbed half-shuffle
fixed-point equations.

Section~\ref{sec:hopf-coalg-structure} establishes the bar-deleting
unshuffle Hopf quotient and describes its kernel in
Theorem~\ref{thm:flattening-morphism} and
Proposition~\ref{prop:hopf-ideal}, respectively. 
It then compares the ordinary-word linear section
$\iota$ with the bar-inserting Hopf algebra section
$\sesec$ of the bar-deleting morphism $\conc$,
as established in
Proposition~\ref{prop:genuine-hopf-section}, 
and studies the associated projections and comodule structures, including the
right-comodule splitting in
Corollary~\ref{cor:comodule-splitting}.

Section~\ref{sec:non-crossing-closure} relates these constructions
to partition combinatorics. After identifying the crossing
correction in~\eqref{eq:local-crossing-correction} and
Lemma~\ref{lem:gap-crossing}, it constructs the classical-to-free
transform and uses Lehner's closure-fiber formula, recalled in
Lemma~\ref{lemma:closure-fiber-formula}, to establish the commutative
diagram of Theorem~\ref{thm:expanded-diagram}.
Example~\ref{ex:degree-four-closure} illustrates the diagram in the
centered degree-four case.
Appendix~\ref{app:low-degree-calculations} gives explicit calculations
of both perturbed recursions in degrees three, four, and five.

\medskip
\noindent\textbf{Notation.}
Throughout the paper, $\bfk$ is a field of characteristic zero,
and $(\A,\cdot_{\A})$ is a nonzero unital associative
$\bfk$-algebra. All vector spaces, tensor products, algebraic
structures, coalgebraic structures, and linear maps are taken over $\bfk$.
For a $\bfk$-vector space $V$, we write $\Lin(V,\bfk)$ for the
space of $\bfk$-linear functionals on $V$.
For a finite set $X$, we write $\bfk^X$ for the space of all
functions $X\longrightarrow\bfk$.
For the coproduct $\Delta$ of a given Hopf algebra $H$, we use $\bar{\Delta}$ to denote the reduced coproduct of $\Delta$. 

\section{Perturbed half-shuffle equations}
\label{sec:perturbed-equations}

In this section, we apply the concatenation perturbation to the moment-perturbed fixed-point equation and to the cumulant-perturbed fixed-point equation. We show that both
constructions give the same classical moment--cumulant relation on
ordinary words.

\subsection{The tensor and double tensor Hopf algebras}
\label{subsec:construction-dcanc}
In this subsection, we recall the tensor and double tensor Hopf algebras, their left and right half-coproducts, and the associated half-shuffle products. Set
\[
\begin{aligned}
 T(\A):=  \bigoplus_{n\geq1}\A^{\otimes n},
 \qquad \overline T(\A):=\bfk\mathbf1\oplus T(\A), \qquad 
 H(\A):= \overline T(T(\A))
 =\bfk\mathbf1\oplus\bigoplus_{s\geq1}T(\A)^{\otimes s}.
\end{aligned}
\]
We write ordinary words as \(a_1\cdots a_n\in \A^{\otimes n}\), with
\(a_1,\ldots,a_n\in\A\), and bar words as \(w_1|\cdots|w_r\in T(\A)^{\otimes s}\),
with \(w_1,\ldots,w_r\in T(\A)\). The tensor algebra $\overline T(\A)$ and the double tensor algebra
$H(\A)$ carry the following connected graded Hopf algebra structures,
respectively~\cite{EFP2015,EFP2019}:
\[
 \bigl(\overline T(\A),m_{\overline T(\A)},\deshat,
       \eta_{\overline T(\A)},\varepsilon_{\overline T(\A)}\bigr),\qquad
 \bigl(H(\A),m_{H(\A)},\deshah,
       \eta_{H(\A)},\varepsilon_{H(\A)}\bigr).
\]
These structures are described explicitly as follows.
\begin{itemize}
\item
The multiplications $m_{\overline T(\A)}$ and $m_{H(\A)}$
are given by ordinary word concatenation and bar concatenation,
respectively. Both algebras are graded by the total number of
letters.

\item
The unit maps
\[
\eta_{\overline T(\A)}:\bfk\longrightarrow\overline T(\A),
\qquad
\eta_{H(\A)}:\bfk\longrightarrow H(\A)
\]
are both defined by $\lambda\mapsto\lambda\mathbf1$.
The corresponding counits
\[
\varepsilon_{\overline T(\A)}:
\overline T(\A)\longrightarrow\bfk,
\qquad
\varepsilon_{H(\A)}:
H(\A)\longrightarrow\bfk
\]
extract the coefficient of the empty word and the empty bar
word, respectively.

\item
To describe the coproducts, let
$I=\{i_1<\cdots<i_n\}\subseteq\mathbb Z_{\geq1}$
be a finite ordered set, where $<$ denotes the usual
strict order on the integers, and let $(a_i)_{i\in I}$
be a family of elements of $\A$. For each subset $S\subseteq I$, write
\begin{equation}\label{eq:concom}
I\setminus S:=J_1\sqcup\cdots\sqcup J_\ell,
\end{equation}
where $J_1,\ldots,J_\ell$ are the maximal consecutive blocks
of $I\setminus S$ with respect to the order induced by $I$.
These blocks are called the \emph{complement gaps} of $S$
in $I$. We use the notation
\begin{equation}\label{eq:aI-conv}
a_I:=a_{i_1}\cdots a_{i_n},
\qquad
a_I^S:=a_{J_1}|\cdots|a_{J_\ell},
\end{equation}
with the conventions $a_\varnothing:=\mathbf1$ and
$a_I^I:=\mathbf1$. The coproducts are then determined on ordinary words by
\begin{equation}\label{eq:coproducts}
\begin{aligned}
\deshat:\quad&
\overline T(\A)\longrightarrow
\overline T(\A)\otimes\overline T(\A),
&
a_I&\longmapsto
\sum_{S\subseteq I}a_S\otimes a_{I\setminus S},
\\
\deshah:\quad&
H(\A)\longrightarrow H(\A)\otimes H(\A),
&
a_I&\longmapsto
\sum_{S\subseteq I}a_S\otimes a_I^S.
\end{aligned}
\end{equation}
Both coproducts map $\mathbf1$ to
$\mathbf1\otimes\mathbf1$ and extend multiplicatively
with respect to their corresponding products, using
componentwise multiplication in the tensor products.
\end{itemize}

Let $H(\A)_n$ and $\overline T(\A)_n$ denote the homogeneous
components of degree $n$ of $H(\A)$ and $\overline T(\A)$,
respectively. We write
\[
H(\A)^+:=\ker\varepsilon_{H(\A)}
=\bigoplus_{n\geq1}H(\A)_n
\]
for the augmentation ideal of $H(\A)$.
The half-coproducts and their dual half-shuffle products,
which play a central role in the double tensor Hopf algebra
framework, are recalled below~\cite{EFP2015,EFP2019}.

\begin{itemize}
\item
For a nonempty finite ordered set $I$, let $i_0=\min I$.
The left and right half-coproducts on ordinary words are
defined, respectively, by
\begin{equation}\label{eq:half-coproducts-H}
 \begin{aligned}
 \DprecH &: T(\A)\longrightarrow T(\A)\otimes H(\A),
 \qquad a_I\longmapsto\sum_{i_0\in S\subseteq I}a_S\otimes a_I^S,\\
 \DsuccH &: T(\A)\longrightarrow\overline T(\A)\otimes H(\A)^+,
 \qquad a_I\longmapsto\sum_{i_0\notin S\subset I}
                         a_S\otimes a_I^S.
 \end{aligned}
\end{equation}

\item
These half-coproducts extend to nonempty bar words by
\begin{equation}\label{eq:left-half-extension}
 \begin{aligned}
 \DprecH : H(\A)^+&\longrightarrow H(\A)^+\otimes H(\A),\\
w_1|\cdots|w_r&\longmapsto
   \DprecH(w_1)|\deshah(w_2)|\cdots|\deshah(w_r),\\[2pt]
 \DsuccH : H(\A)^+&\longrightarrow H(\A)\otimes H(\A)^+,\\
w_1|\cdots|w_r&\longmapsto
   \DsuccH(w_1)|\deshah(w_2)|\cdots|\deshah(w_r).
 \end{aligned}
\end{equation}
In particular, these extensions satisfy
\[
\deshah|_{H(\A)^+}=\DprecH+\DsuccH.
\]

\item  Similarly, with \(a_I^S\) in~\eqref{eq:half-coproducts-H} replaced by
\(a_{I\setminus S}\), there are
\[
\begin{aligned}
\DprecTbar:\quad&
T(\A)\longrightarrow T(\A)\otimes\overline T(\A),
&
a_I&\longmapsto
\sum_{i_0\in S\subseteq I}
a_S\otimes a_{I\setminus S},
\\
\DsuccTbar:\quad&
T(\A)\longrightarrow\overline T(\A)\otimes T(\A),
&
a_I&\longmapsto
\sum_{i_0\notin S\subseteq I}
a_S\otimes a_{I\setminus S}.
\end{aligned}
\]

\item
For linear functionals \(f,g:H(\A)\to\bfk\), the dual half-shuffle
products are the linear maps
\[
\begin{aligned}
 f\prec g &: H(\A)^+\longrightarrow\bfk,
 \qquad x\longmapsto m_{\bfk}(f\otimes g)\DprecH(x),\\
 f\succ g &: H(\A)^+\longrightarrow\bfk,
 \qquad x\longmapsto m_{\bfk}(f\otimes g)\DsuccH(x).
\end{aligned}
\]
where $m_{\bfk}:\bfk\otimes\bfk\longrightarrow\bfk$
denotes multiplication in the ground field $\bfk$.
When at least one of $f$ and $g$ vanishes on $\mathbf1$,
we extend both half-shuffle products to $H(\A)$ by
setting
\[
(f\prec g)(\mathbf1)
=(f\succ g)(\mathbf1)=0.
\]
\end{itemize}

The embedding \(\iota\) is denoted by
\begin{equation}\label{eq:iota-embedding}
 \iota:\overline T(\A)\longrightarrow H(\A).
\end{equation}
We identify ordinary words with their images under \(\iota\) when
the ambient space is clear.

\subsection{Moment and cumulant perturbations}
\label{subsec:left-moment-perturbation}
In this subsection, we introduce the concatenation pertuebation and prove that the perturbed left and cumulant-perturbed fixed-point equations characterize classical cumulants on ordinary words. 

To establish the algebraic framework for the subsequent constructions, let us first introduce three essential ingredients.

\begin{itemize}
\item Define the linear \emph{bar-deleting} map
\begin{equation}\label{eq:conc-map}
 \conc:H(\A)\longrightarrow\overline T(\A),
 \qquad w_1|\cdots|w_r\mapsto w_1\cdots w_r,\quad r\geq 0.
\end{equation}
Using \(\iota\) from~\eqref{eq:iota-embedding}, define the \emph{concatenation perturbation} by
\begin{equation}\label{eq:concatenation-perturbation}
 \ioconc:=\iota\circ\conc-\idH:H(\A)\longrightarrow H(\A),\qquad w_1|\cdots|w_r\mapsto w_1\cdots w_r-w_1|\cdots|w_r, \quad r\geq 0.
\end{equation}
These maps satisfy
\[
 (\idH+\ioconc)(a_I^S)
 =\iota\circ\conc(a_I^S)
 =\iota(a_{I\setminus S}).
\]

\item 
For a linear functional \(c:\overline T(\A)\to\bfk\) with
\(c(\mathbf1)=0\), its \emph{infinitesimal-character extension}
is the linear map \(c_H:H(\A)\to\bfk\) given by
\begin{equation}\label{eq:ext-infcha}
 c_H(\mathbf1):=0, \qquad c_H(w):=c(w), \qquad c_H(w_1|\cdots|w_r):=0,
\end{equation}
where \(w, w_1,\dots,w_r\in T(\A),s\geq 2\).
For \(a_I\) in~\eqref{eq:aI-conv}, write
\begin{equation}\label{eq:mom-cum-func}
  m_I:=\Phi(a_I),\qquad c_I:=c(a_I).
\end{equation}

\item The convolution product on \(\Lin(H(\A),\bfk)\) is
\[
 f*_{H(\A)}g
 :=
 m_{\bfk}\circ(f\otimes g)\circ\deshah.
\]
For a unital linear functional \(\varphi:\A\to\bfk\), let
\(\Phi:H(\A)\to\bfk\) be its \emph{character}~\cite{EFP2019}, defined by
\begin{equation}\label{eq:moment-character}
\Phi(a_1\cdots a_n):=\varphi(a_1\cdot_{\A}\cdots\cdot_{\A}a_n),
 \qquad
 \Phi(w_1|\cdots|w_r):=\prod_{j=1}^r\Phi(w_j).
\end{equation}
Here $n\geq 0$, $r\geq 1$, \( a_1,\ldots,a_n\in\A\), and \(w_1,\ldots,w_r\in T(\A)\).
\end{itemize}

With these preparations, we now use the concatenation
perturbation to define the perturbed moment and cumulant
functionals associated with $\Phi$ and $c_H$, respectively.

\begin{definition}\label{def:perturbed-functionals}
With the concatenation perturbation \(\ioconc:H(\A)\longrightarrow H(\A)\) in~\eqref{eq:concatenation-perturbation}, 
\begin{enumerate} 
\item the \emph{perturbed moment functional} is defined to be
\begin{equation}\label{eq:perturbed-moment-functional}
 \Phiconc:=\Phi+\Phi\circ\ioconc=\Phi\circ\iota\circ\conc:H(\A) \rightarrow \bfk,
\end{equation}
\item the \emph{perturbed cumulant functional} is defined to be
\begin{equation}\label{eq:perturbed-cumulant-functional}
 \cconc:=c_H+c_H\circ\ioconc=c\circ\conc:H(\A) \rightarrow \bfk.
\end{equation}
\end{enumerate}
\end{definition}

Notice that
\[
 \Phiconc(u_1|\cdots|u_r)=\Phi(u_1\cdots u_r),\qquad \cconc(u_1|\cdots|u_r)=c_H(u_1\cdots u_r),\qquad \forall u_1,\ldots,u_r\in T(\A).
\]
\begin{remark}
\begin{enumerate}
\item The functional \(\Phiconc\) need not be a character. Indeed, for nonempty
ordinary words \(u,v\in T(\A)\),
\[
 \Phiconc(u|v)=\Phi(uv),
 \qquad
 \Phiconc(u)\Phiconc(v)=\Phi(u)\Phi(v),\qquad \forall u,v\in T(\A).
\]
Consequently, \(\Phiconc\) is a character precisely when the
restriction of \(\Phi\) to ordinary words is multiplicative for
word concatenation. 
\item The functional \(\cconc\) need not be an infinitesimal character,
since
\[
 \cconc(w_1|\cdots|w_r)=c_H(w_1\cdots w_r),
 \qquad r\geq2,
\]
need not vanish. 
\end{enumerate}
\end{remark}

Following \cite[Section~4]{EFP2019}, we write
\[
G_{H(\A)}(\bfk):=
 \left\{
   \Phi\in\Lin(H(\A),\bfk):
   \begin{array}{l}
     \Phi(\mathbf1)=1,\\
     \Phi(x|y)=\Phi(x)\Phi(y),
     \qquad\forall x,y\in H(\A)
   \end{array}
 \right\}
\]
for the group of characters of \(H(\A)\), with group operation
\(*_{H(\A)}\) and identity \(\varepsilon_{H(\A)}\).
We also write
\[
 g_{H(\A)}(\bfk)
 :=
 \left\{
   \xi\in\Lin(H(\A),\bfk):
   \begin{array}{l}
     \xi(\mathbf1)=0,\\
     \xi(w_1|\cdots|w_r)=0,
     \qquad\forall r\ge2,\ w_1,\ldots,w_r\in T(\A)
   \end{array}
 \right\}
\]
for the Lie algebra of infinitesimal characters of \(H(\A)\),
with bracket
\[
 [\xi,\eta]
 :=
 \xi*_{H(\A)}\eta-\eta*_{H(\A)}\xi,
 \qquad
 \xi,\eta\in g_{H(\A)}(\bfk).
\]
The left half-shuffle logarithm~\cite[Section~4]{EFP2019}
is the map
\[
 \log_{\prec}^{H(\A)}:
  G_{H(\A)}(\bfk)\longrightarrow g_{H(\A)}(\bfk),
 \qquad
 \Phi\longmapsto\log_{\prec}^{H(\A)}(\Phi),
\]
which assigns to each character \(\Phi\in  G_{H(\A)}(\bfk)\)
the unique infinitesimal character
\(\log_{\prec}^{H(\A)}(\Phi)\in g_{H(\A)}(\bfk)\)
satisfying the following fixed-point equation:
\begin{equation}\label{eq:left-half-logarithm}
 \Phi=\varepsilon_{H(\A)}
      +\log_{\prec}^{H(\A)}(\Phi)\prec\Phi.
\end{equation}
To obtain the moment-classical relation in the double tensor Hopf algebra, we define the following perturbation equations.
\begin{definition}
With the perturbed moment and cumulant functionals in Definition~\ref{def:perturbed-functionals}, 
\begin{enumerate} 
  \item the \emph{moment-perturbed (left half-shuffle) fixed-point equation} is defined to be
\begin{equation}\label{eq:classical-left-fixed-point}
 \Phi\circ\iota
 =\varepsilon_{\overline T(\A)}
    +(c_H\prec\Phiconc)\circ\iota
 =\varepsilon_{\overline T(\A)}
    +\bigl(c_H\prec(\Phi+\Phi\circ\ioconc)\bigr)\circ\iota,
\end{equation}
  \item the \emph{cumulant-perturbed (right half-shuffle) fixed-point equation} is defined to be
  \begin{equation}\label{eq:right-fixed-point}
 \Phi\circ\iota
 =\varepsilon_{\overline T(\A)}+(\Phi\succ\cconc)\circ\iota
 =\varepsilon_{\overline T(\A)}
    +\bigl(\Phi\succ(c_H+c_H\circ\ioconc)\bigr)\circ\iota.
\end{equation}
\end{enumerate}
\end{definition}

\begin{remark}
For \(I\neq\varnothing\), let \(i_0:=\min I\). Then, 
\begin{align}
 m_I
 =&\ \Phi(a_I)\notag\\
 =&\ \varepsilon_{H(\A)}(a_I)+(c_H\prec\Phiconc)(a_I)\hspace{1cm}(\text{by }\eqref{eq:classical-left-fixed-point})\notag\\
 =&\ m_{\bfk}(c_H\otimes\Phiconc)\DprecH(a_I)\notag\\
 =&\ \sum_{i_0\in S\subseteq I}c_H(a_S)\Phiconc(a_I^S)\hspace{1cm}(\text{by }\eqref{eq:half-coproducts-H})\notag\\
 =&\ \sum_{i_0\in S\subseteq I}c_S\Phi(a_{I\setminus S})\hspace{1cm}(\text{by }\eqref{eq:perturbed-moment-functional})\notag\\
 =&\ \sum_{i_0\in S\subseteq I}c_Sm_{I\setminus S}.
 \label{eq:classical-first-block-recursion}
\end{align}
Similarly, we also get~\eqref{eq:classical-first-block-recursion} by substituting the right half-coproduct in~\eqref{eq:half-coproducts-H} into~\eqref{eq:right-fixed-point}. The low-degree calculations of \(m_I\) are given in Appendix~\ref{app:left-calculations} and~\ref{app:right-calculations}.
\end{remark}

We now prove that both equations are equivalent
to the classical moment--cumulant formula and determine the
cumulants uniquely. For a finite set \(I\), let \(\Part(I)\) denote
its set of partitions. 

\begin{theorem}\label{thm:classical-perturbations}
Let \(\Phi\) be the character in~\eqref{eq:moment-character} and let
\(c:\overline T(\A)\to \bfk\) be linear with
\(c(\mathbf1)=0\). The following conditions are equivalent.
\begin{enumerate} 
\item\label{item:mom-term-pereq} The moment-perturbed fixed-point equation~\eqref{eq:classical-left-fixed-point} holds. 
\item\label{item:cum-term-pereq} The cumulant-perturbed fixed-point equation~\eqref{eq:right-fixed-point}
holds.
\item\label{item:mom-cumeq} For every nonempty finite set \(I\subseteq \mathbb{Z}_{\geq 1}\) and
every family \((a_i)_{i\in I}\) in \(\A\),
\[
 m_I=\sum_{\pi\in\Part(I)}\prod_{B\in\pi}c_B,
\]
\end{enumerate}
where \(m_I\) and \(c_B\) are given in~\eqref{eq:mom-cum-func}. 
\end{theorem}

Thus, any one of these three equivalent conditions uniquely determines the linear functional $c:\overline T(\A)\rightarrow\bfk$, whose values on ordinary words are precisely the classical cumulants.

\begin{proof}[Proof of Theorem \ref{thm:classical-perturbations}]
The calculations in~\eqref{eq:classical-first-block-recursion} show that
\[
\Item~\ref{item:mom-term-pereq}\Longleftrightarrow\eqref{eq:classical-first-block-recursion}\Longleftrightarrow \Item~\ref{item:cum-term-pereq}.
\]
It remains to prove that \(\Item~\ref{item:mom-cumeq}\Longleftrightarrow\eqref{eq:classical-first-block-recursion}\).   We argue by induction on \(|I|\geq 1\), using \(m_\varnothing=1\). 
For the initial step of \(|I|=1\), \eqref{eq:classical-first-block-recursion} gives \(m_I=c_I\), and so \(\Item ~\ref{item:mom-cumeq}\Longleftrightarrow\eqref{eq:classical-first-block-recursion}\). 
For the inductive step of \(I\neq\varnothing\), put \(i_0:=\min I\). Then the required equivalence follows from 
\begin{align*}
 \sum_{i_0\in S\subseteq I}c_Sm_{I\setminus S}
 =\sum_{i_0\in S\subseteq I}c_S
      \sum_{\sigma\in\Part(I\setminus S)}\prod_{B\in\sigma}c_B
 =\sum_{\pi\in\Part(I)}\prod_{B\in\pi}c_B.
\end{align*}
Here the second equality follows from the induction hypothesis. For the last equality, each partition of \(I\) has a unique block \(S\) containing \(i_0\) and its remaining blocks form a partition of
\(I\setminus S\). 
\end{proof}

We end this section with the following remark. 

\begin{remark}
The restriction to ordinary words is essential for
both equations~\eqref{eq:classical-left-fixed-point} and~\eqref{eq:right-fixed-point}. 
For example, take $a_1a_2|a_3\in H(\A)$ to be a bar word. Then 
it follows from \eqref{eq:left-half-extension} that 
\[
\begin{aligned}
 (c_H\prec\Phiconc)(a_1a_2|a_3)
 =&\ c_{\{1,2\}}m_{\{3\}}+c_{\{1\}}m_{\{2,3\}},\\
 (\varepsilon_{H(\A)} +  c_H\prec\Phiconc)(a_1a_2|a_3)-\Phi(a_1a_2|a_3)
 =&\ c_{\{1\}}c_{\{2,3\}}.
\end{aligned}
\]
Similarly, 
\[
\begin{aligned}
 (\Phi\succ\cconc)(a_1a_2|a_3)
 =&\ c_{\{1,2,3\}}+m_{\{2\}}c_{\{1,3\}}+m_{\{3\}}c_{\{1,2\}}
      +m_{\{2\}}m_{\{3\}}c_{\{1\}},\\
 (\varepsilon_{H(\A)} + \Phi\succ\cconc)(a_1a_2|a_3)-\Phi(a_1a_2|a_3)
 =&\ c_{\{1,2,3\}}+m_{\{2\}}c_{\{1,3\}}.
\end{aligned}
\]
Both differences can be nonzero. 
\end{remark}

\section{Hopf-algebraic properties of the concatenation perturbation}
\label{sec:hopf-coalg-structure}

In this section, we study the map \(\ioconc\) used in both functional
perturbations of Definition~\ref{def:perturbed-functionals} from an algebraic point of view. We first realize bar-deleting as a Hopf quotient, compare the two embeddings into
\(H(\A)\), and identify the complementary projections associated with
\(\iota\circ\conc\) and \(-\ioconc\), and then examine their compatibility
with the coproducts and describe the right comodule structure of
\(\ioconc\).

\subsection{The bar-deleting Hopf quotient}
This subsection is devoted to characterize the Hopf-algebraic properties of the bar-deleting map \(\conc\) in~\eqref{eq:conc-map}. Recall the unshuffle bialgebras 
\begin{align*}
&\bigl(\overline T(\A),m_{\overline T(\A)},\deshat,\reddeshat|_{T(\A)}=\redDprecTbar+\redDsuccTbar,\eta_{\overline T(\A)},\varepsilon_{\overline T(\A)}\bigr),\\
&\bigl(H(\A),m_{H(\A)},\deshah,\reddeshah|_{H(\A)^+}=\redDprecH+\redDsuccH,\eta_{H(\A)},\varepsilon_{H(\A)}\bigr),
\end{align*}
are given in~\cite[Sections~3~and~7]{EFP2015} and~\cite[Sections~3]{EFP2019}, where the left half-unshuffle  \(\redDprecTbar\) and right half-unshuffle  \(\redDsuccTbar\) are defined by
\[
\begin{aligned}
\redDprecTbar:T(\A)&\rightarrow T(\A)\otimes T(\A),\qquad a_I&\mapsto \DprecTbar(a_I)-a_I\otimes\mathbf1,\\
\redDsuccTbar:T(\A)&\rightarrow T(\A)\otimes T(\A),\qquad a_I&\mapsto \DsuccTbar(a_I)-\mathbf1 \otimes a_I,
\end{aligned}
\]
and the left half-unshuffle  \(\redDprecH\) and right half-unshuffle  \(\redDsuccH\) are defined by
\[
\begin{aligned}
\redDprecH:H(\A)^+&\rightarrow H(\A)^+\otimes H(\A)^+,\qquad w_1|\cdots|w_r&\mapsto \DprecH(w_1|\cdots|w_r)-w_1|\cdots|w_r\otimes\mathbf1,\\
\redDsuccH:H(\A)^+&\rightarrow H(\A)^+\otimes H(\A)^+,\qquad w_1|\cdots|w_r&\mapsto \DsuccH(w_1|\cdots|w_r)-\mathbf1 \otimes w_1|\cdots|w_r.
\end{aligned}
\]

\begin{theorem}\label{thm:flattening-morphism}
\begin{enumerate}
  \item\label{item:conc-hopfsurhom} The bar-deleting map 
 \[ 
 \conc : \bigl(H(\A),m_{H(\A)},\deshah,
       \eta_{H(\A)},\varepsilon_{H(\A)}\bigr)\relbar\joinrel\twoheadrightarrow\bigl(\overline T(\A),m_{\overline T(\A)},\deshat,
       \eta_{\overline T(\A)},\varepsilon_{\overline T(\A)}\bigr)
 \]
 is a surjective graded Hopf algebra morphism. 
  \item The bar-deleting map 
 \[ 
 \begin{aligned}
 \conc : &\bigl(H(\A),m_{H(\A)},\deshah,\reddeshah|_{H(\A)^+}=\redDprecH+\redDsuccH,\eta_{H(\A)},\varepsilon_{H(\A)}\bigr)\\
 \relbar\joinrel\twoheadrightarrow&\bigl(\overline T(\A),m_{\overline T(\A)},\deshat,\reddeshat|_{T(\A)}=\redDprecTbar+\redDsuccTbar,\eta_{\overline T(\A)},\varepsilon_{\overline T(\A)}\bigr)
 \end{aligned}
 \]
 is a surjective unshuffle bialgebra morphism. 
\end{enumerate}
\end{theorem}

\begin{proof}
(1) For any two bar words \(x=w_1|\cdots|w_r,y=v_1|\cdots|v_s\in H(\A)\),
\[
 \conc(x|y)=\conc(w_1|\cdots|w_r|v_1|\cdots|v_s)=w_1\cdots w_rv_1\cdots v_s=\conc(x)\conc(y).
\]
Together with \(\conc(\mathbf1)=\mathbf1\), this proves that \(\conc\) is an algebra homophism.
The embedding \(\iota\) in~\eqref{eq:iota-embedding} satisfies
\[
 \conc\circ\iota=\idTbar,
\]
so \(\conc\) is surjective. For any ordinary word \(a_I \in \overline{T}(\A)\), 
\begin{equation}\label{eq:conc-coproduct-word}
 (\conc\otimes\conc)\deshah(\iota(a_I))
 \overset{\eqref{eq:coproducts}}{=}\sum_{S\subseteq I}a_S\otimes\conc(a_I^S)
 \overset{\eqref{eq:conc-map}}{=}\sum_{S\subseteq I}a_S\otimes a_{I\setminus S}
 \overset{\eqref{eq:coproducts}}{=}\deshat(a_I).
\end{equation}

For any bar word \(x=w_1|\cdots|w_r\in H(\A)\), multiplicativity of
both coproducts and of \(\conc\otimes\conc\) yields
\begin{align}
(\conc\otimes\conc)\deshah(x)
=&\ \prod_{j=1}^r(\conc\otimes\conc)\deshah(w_j)\notag\\
=&\ \prod_{j=1}^r\deshat(w_j)\hspace{1cm}(\text{by }\eqref{eq:conc-coproduct-word})\notag\\
=&\ \deshat(w_1\cdots w_r)\notag\\
=&\ \deshat\conc(x).
\label{eq:conc-cop-lett}
\end{align}
The map \(\conc\) also preserves
total word length and the counits. Thus \(\conc\) is a surjective morphism of graded bialgebras and hence of graded Hopf algebras.

(2) By \Item ~\ref{item:conc-hopfsurhom}, we just need to prove that the compatibility of the bar-deleting map and left/right half-unshuffles respectively. For every finite set \(I\subseteq \mathbb{Z}_{\geq 1}\), 
\begin{equation}\label{eq:conc-left-half-word}
  (\conc\otimes\conc)\DprecH(\iota(a_I))
 =\sum_{i_0\in S\subseteq I}
   a_S\otimes a_{I\setminus S}
 =\DprecTbar(a_I).
\end{equation}
Then for any nonempty bar word \(w_1|\cdots|w_r\in H(\A)\),
\[
\begin{aligned}
&(\conc\otimes \conc)\DprecH(w_1|\cdots|w_r)\\
=&\ (\conc\otimes\conc)\Bigl(
 \DprecH(w_1)|\deshah(w_2)|\cdots|\deshah(w_r)
 \Bigr)\hspace{1cm}(\text{by }\eqref{eq:left-half-extension})\\
=&\ \DprecTbar(w_1)\deshat(w_2)\cdots\deshat(w_r)\hspace{1cm}(\text{by }\eqref{eq:conc-cop-lett}\text{ and }\eqref{eq:conc-left-half-word})\\
=&\ \DprecTbar(w_1\cdots w_r)\\
=&\ \DprecTbar \conc(w_1|\cdots|w_r),
\end{aligned}
\]
which shows that 
\[
\begin{aligned}
(\conc\otimes \conc)\redDprecH(w_1|\cdots|w_r)
=&\ (\conc\otimes \conc)\DprecH(w_1|\cdots|w_r)-(\conc\otimes \conc)(w_1|\cdots|w_r\otimes \mathbf{1})\\
=&\ \DprecTbar \conc(w_1|\cdots|w_r)-w_1\cdots w_r\otimes \mathbf{1}\\
=&\ \redDprecTbar \conc(w_1|\cdots|w_r).
\end{aligned}
\]
Restricting instead to \(S\subseteq I\) with \(i_0\notin S\), and using the same multiplicative
extension, proves 
\[
(\conc\otimes \conc)\circ\redDsuccH=\redDsuccTbar\circ\conc.
\]
Thus \(\conc\) is a surjective unshuffle bialgebra morphism. This completes the proof.
\end{proof}

The kernel of a Hopf algebra morphism over a field is a Hopf ideal
\cite[Theorem~4.3.1(b)]{Sweedler1969}; see also
\cite[Chapter~1]{Montgomery1993} and \cite[Section~7.1]{Radford2012}.
Since \(\conc\) preserves total word length, \(\mathcal K:=\ker\conc\) is therefore
a graded Hopf ideal. The following proposition identifies its generators
as a two-sided ideal.

\begin{proposition}\label{prop:hopf-ideal}
The kernel \(\mathcal K\), as a graded Hopf ideal, is generated as a two-sided ideal by
\begin{equation}\label{eq:bar-relations}
 u|v-uv,\qquad u,v\in\overline T(\A).
\end{equation}
\end{proposition}

\begin{proof}
Let
\[
 \mathcal R
 :=\bigl\langle
   u|v-uv:u,v\in\overline T(\A)
 \bigr\rangle
\]
be the two-sided ideal in \(H(\A)\). Then
\[
 \conc\bigl(u|v-uv\bigr)
 =\conc(u)\conc(v)-\conc(uv)
 =uv-uv=0.
\]
Hence \(\mathcal R\subseteq\mathcal K\). For the reverse inclusion, write such a bar word explicitly as
\[
 x:=w_1|\cdots|w_r,
 \qquad r\geq1,\quad w_1,\ldots,w_r\in T(\A).
\]
Claim that
\[
 x-\conc(x)
 =w_1|\cdots|w_r-w_1\cdots w_r
 \in\mathcal R.
\]
We prove this relation by induction on \(r\geq 1\). If \(r=1\), then
\(x=w_1\), and therefore
\[
 x-\conc(x)=w_1-w_1=0\in\mathcal R.
\]
Assume the relation holds for some \(r-1\) with \(r\geq2\).
Put
\[
 x':=w_1|\cdots|\iota(w_{r-1}),
 \qquad
 W:=w_1\cdots w_{r-1}=\conc(x')\in T(\A).
\]
Inserting \(\iota(W)|w_r\) gives
\[
 x-\conc(x)
 =x'|w_r-\iota(Ww_r)
 =\bigl(x'-\iota(\conc(x'))\bigr)|w_r
   +\bigl(\iota(W)|w_r-\iota(Ww_r)\bigr).
\]
The induction hypothesis and the right ideal property imply
\[
 \bigl(x'-\iota(\conc(x'))\bigr)|w_r\in\mathcal R.
\]
The second summand,
\[
 \iota(W)|w_r-\iota(Ww_r),
\]
is a generating relation in~\eqref{eq:bar-relations}, with \(u=W\) and
\(v=w_r\), since \(W,w_r\in T(\A)\subset\overline T(\A)\). Therefore
$
 x-\conc(x)\in\mathcal R.
$
This completes the induction. For an arbitrary \(x\in\mathcal K\), 
\[
 x=x-0=x-\conc(x)\in\mathcal R.
\]
Thus \(\mathcal K\subseteq\mathcal R\).
\end{proof}

\subsection{The embedding and bar-inserting map}
\label{subsec:two-embeddings}
In this subsection, we first clarify the relationship between the embedding $\iota$ and the bar-inserting map $\sesec$, and then establish the Hopf-algebraic properties of the latter.

The embedding \(\iota\) in~\eqref{eq:iota-embedding} satisfies
\(\conc\circ\iota=\idTbar\) by~\eqref{eq:conc-map}.
It preserves neither multiplication nor the coproduct. For example, for nonzero
\(a,b,c\in\A\), 
\begin{equation}\label{eq:ito-noalg-nocoalg}
\begin{aligned}
 \iota(a)|\iota(b)-\iota(ab)
 =&\ a|b-ab\neq0,\\
 \deshah\iota(abc)-(\iota\otimes\iota)\deshat(abc)
 =&\ b\otimes(a|c-ac)\neq0.
\end{aligned}
\end{equation}
Define the \emph{bar-inserting map}
\[
 \sesec:\overline T(\A)\lhook\joinrel\longrightarrow H(\A),
 \qquad a_{i_1}\cdots a_{i_n}\longmapsto a_{i_1}|\cdots|a_{i_n}.
\]
Note that
\begin{equation}\label{eq:con-ses}
  \conc\circ \sesec=\conc\circ\iota=\idTbar,
 \qquad
 \iota\circ\conc\neq \sesec\circ\conc.
\end{equation}
\begin{proposition}\label{prop:genuine-hopf-section}
\begin{enumerate}
  \item\label{them:sec-hopfinjhom} The bar-inserting map
 \[ 
 \sesec : \bigl(\overline T(\A),m_{\overline T(\A)},\deshat,
       \eta_{\overline T(\A)},\varepsilon_{\overline T(\A)}\bigr)\lhook\joinrel\longrightarrow\bigl(H(\A),m_{H(\A)},\deshah,
       \eta_{H(\A)},\varepsilon_{H(\A)}\bigr)
 \]
 is an injective graded Hopf algebra morphism. 
  \item The bar-inserting map
 \[ 
 \begin{aligned}
 \sesec : &\bigl(\overline T(\A),m_{\overline T(\A)},\deshat,\reddeshat|_{T(\A)}=\redDprecTbar+\redDsuccTbar,\eta_{\overline T(\A)},\varepsilon_{\overline T(\A)}\bigr)\\
 \lhook\joinrel\longrightarrow&\bigl(H(\A),m_{H(\A)},\deshah,
 \reddeshah|_{H(\A)^+}=\redDprecH+\redDsuccH,\eta_{H(\A)},\varepsilon_{H(\A)}\bigr)
 \end{aligned}
 \]
 is an injective unshuffle bialgebra morphism. 
\end{enumerate}
\end{proposition}

\begin{proof}
(1) Let \(u=a_I\) and \(v=b_J\) be any two ordinary words, with
\(I=\{i_1<\cdots<i_n\}\) and \(J=\{j_1<\cdots<j_m\}\). Then
\[
 \sesec(uv)
 =a_{i_1}|\cdots|a_{i_n}|b_{j_1}|\cdots|b_{j_m}
 =\sesec(u)|\sesec(v).
\]
Together with \(\sesec(\mathbf1)=\mathbf1\), this proves \(\sesec\) is an algebra morphism.
Since \(\conc\circ\sesec=\idTbar\), the map \(\sesec\) is injective. Further, 
\[
\begin{aligned}
 \deshah(\sesec(a_I))
 =&\ \deshah(a_{i_1}|\cdots|a_{i_n})\\
 =&\ \prod_{q=1}^n
   \bigl(a_{i_q}\otimes\mathbf1+\mathbf1\otimes a_{i_q}\bigr)\hspace{1cm}(\text{by }\eqref{eq:coproducts})\\
 =&\ \sum_{S\subseteq I}
   \sesec(a_S)\otimes\sesec(a_{I\setminus S})\\
 =&\ (\sesec\otimes\sesec)\deshat(a_I).
\end{aligned}
\]
The map \(\sesec\) preserves total word length and the counits. Thus \(\sesec\) is an
injective morphism of graded bialgebras and hence of graded Hopf algebras.

(2) By \Item ~\ref{them:sec-hopfinjhom}, it remains to verify that \(\sesec\)
is compatible with the left/right half-unshuffles respectively. Indeed, for every finite set \(I\subseteq \mathbb{Z}_{\geq 1}\), 
\[
\begin{aligned}
 \DprecH\sesec(a_I)
 =&\ \DprecH(a_{i_1}|\cdots|a_{i_n})\\
 =&\ \bigl(a_{i_0}\otimes\mathbf1\bigr)
   \prod_{q=2}^n
   \bigl(a_{i_q}\otimes\mathbf1+\mathbf1\otimes a_{i_q}\bigr)\hspace{1cm}(\text{by }\eqref{eq:half-coproducts-H}\text{ and }\eqref{eq:left-half-extension})\\
 =&\ \sum_{i_0\in S\subseteq I}
   \sesec(a_S)\otimes\sesec(a_{I\setminus S})\\
 =&\ (\sesec\otimes\sesec)\DprecTbar(a_I),
\end{aligned}
\]
which implies that
\[
\begin{aligned}
\redDprecH \sesec(a_I)
=&\ \DprecH \sesec(a_I) - a_{i_1}|\cdots|a_{i_n}\otimes \mathbf{1}\\
=&\ (\sesec\otimes\sesec)\DprecTbar(a_I) - (\sesec\otimes\sesec)(a_I\otimes \mathbf{1})\\
=&\ (\sesec\otimes\sesec)\redDprecTbar(a_I).
\end{aligned}
\]
Working now on \(S\subseteq I\) with \(i_0\notin S\), the same
calculation gives
\[
 \redDsuccH\sesec(a_I)=(\sesec\otimes\sesec)\redDsuccTbar(a_I).
\]
Thus \(\sesec\) is an injective unshuffle bialgebra morphism.
This completes the proof.
\end{proof}

The bar-inserting map \(\sesec\) determines an idempotent Hopf endomorphism of \(H(\A)\).
\begin{corollary}\label{cor:Q-hopf}
\begin{enumerate}
  \item\label{item:seccon-hopfidehom} The map 
 \[ 
 \sesec\circ\conc : \bigl(H(\A),m_{H(\A)},\deshah,
       \eta_{H(\A)},\varepsilon_{H(\A)}\bigr)\longrightarrow\bigl(H(\A),m_{H(\A)},\deshah,
       \eta_{H(\A)},\varepsilon_{H(\A)}\bigr)
 \]
 is an idempotent graded Hopf algebra endomorphism. 
  \item The map
 \[ 
 \begin{aligned}
 \sesec\circ\conc : &\bigl(H(\A),m_{H(\A)},\deshah,
 \reddeshah|_{H(\A)^+}=\redDprecH+\redDsuccH,\eta_{H(\A)},\varepsilon_{H(\A)}\bigr)\\
\longrightarrow&\bigl(H(\A),m_{H(\A)},\deshah,
\reddeshah|_{H(\A)^+}=\redDprecH+\redDsuccH,\eta_{H(\A)},\varepsilon_{H(\A)}\bigr)
 \end{aligned}
 \]
 is an idempotent unshuffle bialgebra endomorphism.
\end{enumerate}
\end{corollary}

\begin{proof}
(1) By Theorem~\ref{thm:flattening-morphism} and
Proposition~\ref{prop:genuine-hopf-section}, the composition
\(\sesec\circ\conc\) is a graded Hopf algebra endomorphism of \(H(\A)\). Further, 
\[
 (\sesec\circ\conc)^2
 =\sesec\circ(\conc\circ\sesec)\circ\conc
 \overset{\eqref{eq:con-ses}}{=}\sesec\circ\conc.
\]

(2) It follows from \Item ~\ref{item:seccon-hopfidehom} and
\[
\begin{aligned}
 \redDprecH\circ\sesec\circ\conc
 =&\ (\sesec\otimes\sesec)\circ\redDprecTbar\circ\conc
 \hspace{1cm}(\text{by Proposition }\ref{prop:genuine-hopf-section})\\
 =&\ (\sesec\otimes\sesec)\circ(\conc\otimes\conc)\circ\redDprecH
 \hspace{1cm}(\text{by Theorem }\ref{thm:flattening-morphism})\\
 =&\ \bigl((\sesec\circ\conc)\otimes(\sesec\circ\conc)\bigr)\circ\redDprecH.
\end{aligned}
\]
And in the same way,
\[
 \redDsuccH\circ\sesec\circ\conc
 =\bigl((\sesec\circ\conc)\otimes(\sesec\circ\conc)\bigr)\circ\redDsuccH.
\]
This completes the proof.
\end{proof}

A simple result of Corollary~\ref{cor:Q-hopf} is 
\[
H(\A)=\operatorname{im}(\sesec\circ\conc)\oplus\ker(\sesec\circ\conc)=\sesec(\overline T(\A))\oplus_{g}\mathcal K,
\]
since \(\sesec\circ\conc:H(\A)\rightarrow H(\A)\) is an idempotant operator. We now collect some elementary properties of the endomorphism \(\iota\circ\conc\) and the perturbation
\(\ioconc\) in~\eqref{eq:concatenation-perturbation} via a remark. 
\begin{remark}\label{remark:itocon}
\begin{enumerate}
  \item The map $$\iota\circ\conc:H(\A)\longrightarrow H(\A) $$ is a graded linear idempotent operator, and the map \(\ioconc:H(\A)\longrightarrow H(\A)\) is a graded linear negative idempotent operator, with 
\[
 \operatorname{im}(\iota\circ\conc)=\iota(\overline T(\A)),
 \qquad
 \ker(\iota\circ\conc)=\mathcal K,
 \qquad
 \operatorname{im}\ioconc=\mathcal K,
 \qquad
 \ker\ioconc=\iota(\overline T(\A)).
\]
  \item\label{item:iotcon-proj-stru} We have
\begin{equation}\label{eq:H-splitting}
 H(\A)=\operatorname{im}(\iota\circ\conc)\oplus\ker(\iota\circ\conc)=\iota(\overline T(\A))\oplus_g\mathcal K,
\end{equation}
since \(\iota\circ\conc:H(\A)\longrightarrow H(\A)\) is an idempotent operator. 
\item The decomposition~\eqref{eq:H-splitting} is as graded vector
spaces. The calculations in~\eqref{eq:ito-noalg-nocoalg}
show that \(\iota(\overline T(\A))\) is neither a subalgebra nor a
subcoalgebra of \(H(\A)\). Consequently, the projector \(\iota\circ\conc\) does not
preserve the Hopf structure. The Hopf projector associated with \(\conc\) is
\(\sesec\circ\conc\) from Corollary~\ref{cor:Q-hopf}.
  \item Neither \(\iota\circ\conc\) nor \(\ioconc\) is an algebra morphism
for the bar product on \(H(\A)\). Indeed, for nonzero \(a,b\in\A\),
we have \(\ioconc(a)=\ioconc(b)=0\) and
\[
\begin{aligned}
 (\iota\circ\conc)(a|b)
 =&\ ab\neq a|b=(\iota\circ\conc)(a)|(\iota\circ\conc)(b),\\
 \ioconc(a|b)
 =&\ ab-a|b\neq0=\ioconc(a)|\ioconc(b).
\end{aligned}
\]
  \item The map \(\ioconc\) is not a derivation either, since for the same
\(a,b\), the Leibniz rule fails:
\[
 \ioconc(a|b)=ab-a|b\neq0=\ioconc(a)|b+a|\ioconc(b).
\]
\end{enumerate}
\end{remark}

\subsection{The comodule structure on \(H(\A)\) }
\label{subsec:comodule-structure}
In this subsection, we study the right and left \(\overline T(\A)\)-comodule structures induces by the bar-deleting map \(\conc\). Define
\begin{equation}\label{eq:right-quotient-coaction}
\begin{aligned}
 \rhoR:=&\ (\idH\otimes\conc)\circ\deshah:
 H(\A)\longrightarrow H(\A)\otimes\overline T(\A),\\
 \rhoL:=&\ (\conc\otimes\idH)\circ\deshah:
 H(\A)\longrightarrow\overline T(\A)\otimes H(\A).
\end{aligned}
\end{equation}

\begin{theorem}\label{thm:right/left-coaction}
\begin{enumerate}
  \item\label{item:rightcoa} The pair \((H(\A),\rhoR)\) is a right \(\overline T(\A)\)-comodule algebra.
  \item The pair \((H(\A), \rhoL)\) is a left \(\overline T(\A)\)-comodule algebra. 
\end{enumerate}
\end{theorem}

\begin{proof}
(1) We have
\[
\begin{aligned}
 &(\rhoR\otimes\idTbar)\circ\rhoR\\
 =&\ (\idH\otimes \conc\otimes \conc)\circ(\deshah\otimes\idH)\circ\deshah\\
 =&\ (\idH\otimes \conc\otimes \conc)\circ(\idH\otimes\deshah)\circ\deshah\hspace{1cm}(\text{by coassociativity of }\deshah)\\
 =&\ (\idH\otimes\deshat\circ\conc)\circ\deshah\hspace{4cm}(\text{by }\eqref{eq:conc-cop-lett})\\
 =&\ (\idH\otimes\deshat)\circ(\idH\otimes\conc)\circ\deshah\\
 =&\ (\idH\otimes\deshat)\circ\rhoR.
\end{aligned}
\]
Also,
\[
 (\idH\otimes\varepsilon_{\overline T(\A)})\circ\rhoR
 =(\idH\otimes\varepsilon_{\overline T(\A)}\conc)\circ\deshah
 =(\idH\otimes\varepsilon_{H(\A)})\circ\deshah
 =\idH.
\]
Since \(\deshah\) and \(\idH\otimes\conc\) are algebra morphisms, their
composition \(\rhoR\) is also an algebra morphism. 

(2) It is similar to the proof of \Item ~\ref{item:rightcoa}. 
\end{proof}

We equip \(\overline T(\A)\) with its regular right coaction \(\deshat\). The following proposition shows that \(\iota\) and \(\conc\) are right \(\overline T(\A)\)-comodule morphisms.
\begin{proposition}\label{prop:right-colinearity}
The maps
  \[
  \iota:(\overline T(\A),\deshat)\hookrightarrow(H(\A),\rhoR),
  \qquad
  \conc:(H(\A),\rhoR)\twoheadrightarrow(\overline T(\A),\deshat),
  \]
are right \(\overline T(\A)\)-comodule morphisms.
\end{proposition}

\begin{proof}
Since \(\conc\circ\iota=\idTbar\), the map \(\iota\) is injective
and \(\conc\) is surjective. It remains to verify their right-colinearity.
For any nonempty ordinary word \(a_I\in T(\A)\), 
\[
\begin{aligned}
 (\rhoR\circ\iota)(a_I)
 =&\ (\idH\otimes\conc)\circ\deshah\circ(\iota(a_I))
 \hspace{1cm}(\text{by }\eqref{eq:right-quotient-coaction})\\
 =&\ \sum_{S\subseteq I}\iota(a_S)\otimes\conc(a_I^S)
 \hspace{1cm}(\text{by }\eqref{eq:coproducts})\\
 =&\ \sum_{S\subseteq I}\iota(a_S)\otimes a_{I\setminus S}
 \hspace{1cm}(\text{by }\eqref{eq:conc-map}\text{ and }\eqref{eq:aI-conv})\\
 =&\ (\iota\otimes\idTbar)\circ\deshat(a_I).
 \hspace{1cm}(\text{by }\eqref{eq:coproducts})
\end{aligned}
\]
The same identity holds on the unit, so by linearity,
\[
 \rhoR\circ\iota=(\iota\otimes\idTbar)\circ\deshat.
\]
Thus \(\iota\) is a right \(\overline T(\A)\)-comodule morphism.
For \(\conc\), 
\[
\begin{aligned}
 (\conc\otimes\idTbar)\circ\rhoR
 =&\ (\conc\otimes\idTbar)\circ(\idH\otimes\conc)\circ\deshah
 \hspace{1cm}(\text{by }\eqref{eq:right-quotient-coaction})\\
 =&\ (\conc\otimes\conc)\circ\deshah\\
 =&\ \deshat\circ\conc
 \hspace{1cm}(\text{by }\eqref{eq:conc-cop-lett}).
\end{aligned}
\]
Hence \(\conc\) is also a right \(\overline T(\A)\)-comodule morphism. This completes the proof.
\end{proof}
\begin{remark}
The maps
  \[
  \iota\circ\conc:(H(\A),\rhoR)\longrightarrow(H(\A),\rhoR),
  \qquad
  \ioconc:(H(\A),\rhoR)\longrightarrow(H(\A),\rhoR)
  \]
are endomorphisms of the right \(\overline T(\A)\)-comodule \((H(\A),\rhoR)\).
\end{remark}

The direct sum decomposition from Remark~\ref{remark:itocon} therefore split \(H(\A)\) as a
right comodule. 
\begin{corollary}\label{cor:comodule-splitting}
We have the following split short exact sequence of right \(\overline T(\A)\)-comodules
\[
\begin{tikzcd}[column sep=3.2em]
 0 \arrow[r]
 & \mathcal K \arrow[r,hook,"\jmath_{\mathcal K}"]
 & H(\A) \arrow[r,shift left=.7ex,"\conc"]
 & \overline T(\A)
   \arrow[l,shift left=.7ex,"\iota"]
   \arrow[r]
 & 0.
\end{tikzcd}
\]
\end{corollary}

\begin{proof}
By Remark~\ref{remark:itocon}, \(\mathcal K=\operatorname{im}\ioconc\). Since \(\ioconc\) is a right \(\overline T(\A)\)-comodule morphism by Proposition~\ref{prop:right-colinearity}, \(\mathcal K\) is a right subcomodule of \((H(\A),\rhoR)\). Since \(\jmath_{\mathcal K}\) is injective, \(\operatorname{im}\jmath_{\mathcal K}
=\mathcal K=\ker\conc\), and \(\conc\) is surjective, the sequence is a short exact sequence of right \(\overline T(\A)\)-comodules. Finally, \(\conc\circ\iota=\idTbar\), as required.
\end{proof}

\begin{remark}
\begin{enumerate}
  \item Nevertheless, the embedding \(\iota\), the projector \(\iota\circ\conc\), and
the concatenation perturbation \(\ioconc\) are, in general, not left
\(\overline T(\A)\)-comodule morphism.
For example, choose
nonzero \(a,b,c\in\A\), and set \(x=abc\). Then
\[
\begin{aligned}
 &\rhoL\circ\iota(abc)-(\idTbar\otimes\iota)\circ\deshat(abc)
 =b\otimes(a|c-ac)\neq0,\\
 &\rhoL\circ(\iota\circ\conc)(x)-(\idTbar\otimes(\iota\circ\conc))\circ\rhoL(x)
 =b\otimes(a|c-ac)\neq0,\\
 &(\idTbar\otimes\ioconc)\circ\rhoL(x)
 =b\otimes(ac-a|c)\neq0=\rhoL\circ\ioconc(x).
\end{aligned}
\]
  \item\label{item:bar-ins-rigcom} The bar-inserting map
  \[
   \sesec:(\overline T(\A),\deshat)\hookrightarrow(H(\A),\rhoR)
  \]
  is a right \(\overline T(\A)\)-comodule morphism since
  \[
\begin{aligned}
 \rhoR\circ\sesec
 =&\ (\idH\otimes\conc)\circ\deshah\circ\sesec
 \hspace{1cm}(\text{by }\eqref{eq:right-quotient-coaction})\\
 =&\ (\idH\otimes\conc)\circ(\sesec\otimes\sesec)\circ\deshat
 \hspace{1cm}(\text{by Proposition }\ref{prop:genuine-hopf-section})\\
 =&\ (\sesec\otimes\conc\circ\sesec)\circ\deshat\\
 =&\ (\sesec\otimes\idTbar)\circ\deshat.
 \hspace{1cm}(\text{by }\eqref{eq:con-ses})
\end{aligned}
\]
  \item The bar-inserting
  \[
   \sesec:(\overline T(\A),\deshat)\hookrightarrow(H(\A),\rhoL)
  \]
  is a left \(\overline T(\A)\)-comodule morphism, which is similar to \Item ~\ref{item:bar-ins-rigcom}.
\end{enumerate}
\end{remark}

\section{Classical-to-free cumulant transform}
\label{sec:non-crossing-closure}
In this section, we relate the classical-to-free cumulant transform
to push-forward along non-crossing closure through functions
on partitions,  using the closure-fiber formula in
\cite[Theorem~4.1 and its proof]{Lehner2002}.

\subsection{Non-crossing closure and the classical-to-free cumulant transform}
\label{subsec:classical-to-free-transform}
Let us recall from~\cite{Lehner2002} how free cumulants can be expressed in terms of classical cumulants via non-crossing closure and connected partitions. 

Fix a finite \(I\subseteq \mathbb{Z}_{\geq 1}\), and recall \(\Part(I)\) is the set of partitions of \(I\). Distinct blocks \(B,C\) cross if there are \(b_1<b_2\) in \(B\)
and \(c_1<c_2\) in \(C\) with \(b_1<c_1<b_2<c_2\) or
\(c_1<b_1<c_2<b_2\). Let \(\NC(I)\) denote the partitions of
\(I\) with no crossing pair of blocks. For partitions, \(\sigma\leq\pi\)
means that every block of \(\sigma\) is contained in a block of \(\pi\).

\begin{definition}[{\cite[Theorem~4.1 and its proof]{Lehner2002}}]
\label{def:closure-pushforward}
The \emph{non-crossing closure map} is
\[
 \clNC^I:\Part(I)\longrightarrow\NC(I),
 \qquad
 \sigma\longmapsto
 \overline\sigma:=\min\{\pi\in\NC(I):\sigma\leq\pi\}.
\]
\end{definition}

\noindent The \emph{closure push-forward} is the \(\bfk\)-linear map
\[
 (\clNC^I)_*:\bfk^{\Part(I)}\longrightarrow\bfk^{\NC(I)},
 \qquad w\longmapsto(\clNC^I)_*w,
\]
defined by the standard summation over fibers
\cite[Definition~2.5]{BaezCourser2018}:
\[
 \bigl((\clNC^I)_*w\bigr)(\pi)
 :=\sum_{\substack{\sigma\in\Part(I)\\
                   \clNC^I(\sigma)=\pi}}w(\sigma),
 \qquad \pi\in\NC(I).
\]

The superscript \(I\) is omitted from both maps when the ordered set is clear. For \(I\neq\varnothing\), put \(\onepart{I}:=\{I\}\), and let
\(\Part^{\Conn}(I)\) consist of the \emph{connected partitions}:
no nonempty proper interval of \(I\) is a union of their blocks.
Equivalently, these are precisely the partitions
\(\sigma\in\Part(I)\) such that
\(\overline\sigma=\onepart{I}\);
see \cite[proof of Theorem~4.1]{Lehner2002}.

For a linear functional \(f:\overline T(\A)\to\bfk\), write
\begin{equation}\label{eq:block-products}
 f_B:=f(a_B),
 \qquad f_\pi:=\prod_{B\in\pi}f_B,
 \qquad \pi\in\Part(I),
\end{equation}
with the empty product equal to \(1\). The following remark gives a factorization
of the sum over each closure fiber.

\begin{remark}[{\cite[proof of Theorem~4.1]{Lehner2002}}]
\label{remark:fiber-factorization}
For every \(\pi\in\NC(I)\), the products
\(c_\sigma\) defined in~\eqref{eq:block-products} satisfy
\[
 \sum_{\substack{\sigma\in\Part(I)\\
                  \overline\sigma=\pi}}c_\sigma
 =
 \prod_{B\in\pi}
 \left(
   \sum_{\tau\in\Part^{\Conn}(B)}c_\tau
 \right).
\]
\end{remark}
We now construct the classical-to-free cumulant transform. Let \(c:\overline T(\A)\to\bfk\) be any linear functional with
\(c(\mathbf1)=0\), and let \(c_H\) be its infinitesimal-character
extension to \(H(\A)\), as in~\eqref{eq:ext-infcha}. Define
\(\Phi_c:H(\A)\to\bfk\) as the bar-multiplicative extension of the
classical word moments
\begin{equation}\label{eq:classical-word-moments}
 \Phi_c(a_I):=\sum_{\sigma\in\Part(I)}c_\sigma,
\end{equation}
where \(c_\sigma\) are given in~\eqref{eq:block-products}.
Explicitly, 
\[
 \Phi_c(w_1|\cdots|w_r):=\prod_{j=1}^r\Phi_c(w_j),\quad \text{with the convention that }\Phi_c(\mathbf1) :=1.
\]
Thus \(\Phi_c\) is a character. When \(c\) is the classical cumulant word functional of \(\varphi\) in~\eqref{eq:moment-character}, the classical moment--cumulant formula and bar-multiplicativity give \(\Phi_c=\Phi\). Throughout this section, write
\(m_I:=\Phi_c(a_I)\). 

We interpret moments and classical cumulants within the framework of double tensor Hopf algebra as follows. 
\begin{itemize}
\item The tensor-algebra convolution product on $\Lin(\overline T(\A),\bfk)$ and its exponential are~\cite[Sections~7~and~8]{EFP2015}
\begin{equation}\label{eq:tens-con}
  f\convTbar g:=m_{\bfk}(f\otimes g)\deshat,
 \qquad
 \mu_c:=\exp^{\convTbar}(c)
 =\varepsilon_{\overline T(\A)}
   +\sum_{q\geq1}\frac{1}{q!}c^{\convTbar q}.
\end{equation}

\item For a unital linear functional \(\mu:\overline T(\A)\rightarrow \bfk\), let
\(\Ext(\mu)\in G_{H(\A)}(\bfk)\) be its bar-multiplicative extension:
\[
 \Ext(\mu)(w_1|\cdots|w_r):=\prod_{j=1}^r\mu(w_j).
\]

\item Define the spaces
\[
 \begin{aligned}
 \Lin_0(\overline T(\A),\bfk):=&\
 \{c\in \Lin(\overline T(\A),\bfk):c(\mathbf1)=0\}.
 \end{aligned}
\]
For $c\in  \Lin_0(\overline T(\A),\bfk) $, 
\begin{equation}\label{eq:muc-phic}
 \mu_c\overset{\eqref{eq:classical-word-moments}}{=}\Phi_c\circ\iota:\overline T(\A)\rightarrow\bfk,\qquad \Phi_c\overset{\eqref{eq:classical-word-moments}}{=}\Ext(\mu_c)\in G_{H(\A)}(\bfk).
\end{equation}

\item For $c\in  \Lin_0(\overline T(\A),\bfk) $, by Definition~\ref{def:perturbed-functionals},  
\begin{equation}\label{eq:phic-con}
  \Phi_c^{\mathfrak{c}}:=\Phi_c+\Phi_c\circ\ioconc=\Phi_c\circ\iota\circ\conc=\mu_c\circ\conc.
\end{equation}
Replacing $\Phi$ with $\Phi_c$ in the moment-perturbed fixed-point equation~\eqref{eq:classical-left-fixed-point} and the
cumulant-perturbed fixed-point equation~\eqref{eq:right-fixed-point}, we obtain
\begin{equation}\label{eq:classical-word-recursion}
 \Phi_c\circ\iota
 =\varepsilon_{\overline T(\A)}
       +(c_H\prec\Phi_c^{\mathfrak{c}})\circ\iota
 =\varepsilon_{\overline T(\A)}
       +(\Phi_c\succ\cconc)\circ\iota.
\end{equation}
\end{itemize}

\begin{remark}
We now identify the crossing terms contributed by \(\ioconc\)
to the moment-perturbed identity in~\eqref{eq:classical-word-recursion}. Fix a nonempty finite set \(I\subseteq \mathbb{Z}_{\geq 1}\) and a proper subset
\(S\subsetneq I\) with \(\min I\in S\). Under~\eqref{eq:concom}~and~\eqref{eq:aI-conv}, the corresponding term of the left half-coproduct in~\eqref{eq:half-coproducts-H} has second tensor factor
\(a_{J_1}|\cdots|a_{J_\ell}\). Then 
\begin{align}
 &\ (\Phi_c\circ\ioconc)(a_{J_1}|\cdots|a_{J_\ell})\notag\\
 =&\ \Phi_c^{\mathfrak{c}}(a_{J_1}|\cdots|a_{J_\ell})-\Phi_c(a_{J_1}|\cdots|a_{J_\ell})\notag\\
 =&\ m_{I\setminus S}-\prod_{q=1}^\ell m_{J_q}\notag\\
 =&\ \sum_{\tau\in\Part(I\setminus S)}c_\tau
 -\sum_{\substack{\tau\in\Part(I\setminus S)\\
                   \forall B\in\tau\ \exists q:\ B\subseteq J_q}}
   c_\tau\notag\\
 =&\ \sum_{\substack{\tau\in\Part(I\setminus S)\\
        \exists B\in\tau,\ \exists p\neq q:\\
        B\cap J_p\neq\varnothing,\ B\cap J_q\neq\varnothing}}
   c_\tau.
\label{eq:local-crossing-correction}
\end{align}
Thus \(\Phi_c\circ\ioconc\) retains exactly the terms \(c_\tau\) for
complement partitions with a block joining two different gaps. The following
elementary consequence of the standard first-block decomposition
of non-crossing partitions identifies this condition with crossing
\(S\).
\end{remark}

\begin{lemma}[{\cite[proof of Theorem~16.15]{NS2006}
and \cite[proof of Theorem~5.2]{EFP2015}}]
\label{lem:gap-crossing}
Let \(I\) be a nonempty finite linearly ordered set, and let
\(S\subseteq I\) contain \(\min I\).
Let \(J_1,\ldots,J_\ell\) be the connected components of
\(I\setminus S\) relative to \(I\).
For every nonempty subset \(B\subseteq I\setminus S\),
the following conditions are equivalent:
\begin{enumerate}
\item \(B\) meets two distinct complement gaps; that is,
\(B\cap J_p\neq\varnothing\) and \(B\cap J_q\neq\varnothing\)
for some \(p\neq q\).
\item The blocks \(B\) and \(S\) cross.
\end{enumerate}
\end{lemma}

Combining Proposition~\ref{lem:gap-crossing} with
\eqref{eq:local-crossing-correction}, and multiplying by the
factor \(c_S\), gives
\[
 c_S\cdot_{\bfk}(\Phi_c\circ\ioconc)(a_{J_1}|\cdots|a_{J_\ell})
 =\sum_{\substack{\tau\in\Part(I\setminus S)\\
                  \exists B\in\tau:\ B\text{ crosses }S}}
   c_{\{S\}\sqcup\tau}.
\]
For the fixed block \(S\), the terms produced by the perturbation \(\ioconc\) is the sum of the terms corresponding to partitions of \(I\) that contain \(S\) and another block crossing \(S\). This sum is zero if \(I\setminus S\) has at most one connected component. Applying the same recursion to the moments of the complement also includes crossings between the remaining blocks. 

We next apply the left half-shuffle logarithm to \(\Phi_c\) to obtain its free cumulants $r_c$. 

\begin{definition}
The \emph{classical-to-free transform} is defined as
\begin{equation}\label{eq:classical-to-free-transform}
 \mathcal L : \Lin_0(\overline T(\A),\bfk) \rightarrow g_{H(\A)}(\bfk),\qquad c\mapsto r_c:=\log_{\prec}^{H(\A)}(\Phi_c).
\end{equation}
\end{definition}

On ordinary words, the half-shuffle expansion
\cite[Theorem~5.2]{EFP2015} yields Speicher's non-crossing
moment--cumulant formula \cite{Speicher1994}:
\[
 m_I=\Phi_c(a_I)=\sum_{\pi\in\NC(I)}(r_c)_\pi.
\]

We recall the following multilinear form of Lehner's closure-fiber
formula. Its multivariate version is recorded in
\cite[Section~1]{AHLV2015}. In the present notation, it expresses
the products \((r_c)_\pi\) associated with the transform
\eqref{eq:classical-to-free-transform} as sums of the products
\(c_\sigma\) over closure fibers.

\begin{lemma}[{\cite[Theorem~4.1 and its proof]{Lehner2002}}]
\label{lemma:closure-fiber-formula}
Let \(c\in\Lin_0(\overline T(\A),\bfk)\) and \(r_c=\mathcal L(c)\).
Then, for every finite set \(I\subseteq \mathbb{Z}_{\geq 1}\), every
family \((a_i)_{i\in I}\) of elements of \(\A\), and every \(\pi\in\NC(I)\),
\[
 (r_c)_\pi=\sum_{\substack{\sigma\in\Part(I)\\
                  \overline\sigma=\pi}}c_\sigma.
\]
In particular, when \(\pi= \{I\}\),
\[
(r_c)_\pi = r_c(a_I)=\sum_{\sigma\in\Part^{\Conn}(I)}c_\sigma.
\]
\end{lemma}
\subsection{Compatibility with non-crossing closure}
In this subsection, we use a commutative diagram to show how the bar-deleting map \(\conc\) in~\eqref{eq:conc-map} and the classical-to-free cumulant transform \(\mathcal L\) in~\eqref{eq:classical-to-free-transform} recover Lehner's closure-fiber formula recalled in Lemma~\ref{lemma:closure-fiber-formula}. 

\begin{definition}\label{def:partition-function-maps}
Fix a finite set \(I\subseteq \mathbb{Z}_{\geq 1}\) and a family
\((a_i)_{i\in I}\) in \(\A\). Define the maps 
\[
P^{\mathrm{cla}}:\Lin_0(\overline T(\A),\bfk)\longrightarrow \bfk^{\Part(I)},
 \qquad c\mapsto P^{\mathrm {cla}}(c):= c_{\mybullet} : \Part(I)\rightarrow \bfk, \qquad \sigma\mapsto c_\sigma,
\]
and
\[
 P^{\mathrm{fre}}: g_{H(\A)}(\bfk)\longrightarrow \bfk^{\NC(I)},
 \qquad r\mapsto P^{\mathrm {fre}}(r):= r_{\mybullet} : \NC(I)\rightarrow \bfk, \qquad \pi\mapsto r_\pi,
\]
Here the products \( c_\sigma\) and \(r_\pi\) are defined in~\eqref{eq:block-products}. 
\end{definition}

To display the concatenation perturbation explicitly in the diagram, put
\[
 \Lin_0(H(\A),\bfk)
 :=\{\beta\in\Lin(H(\A),\bfk):\beta(\mathbf1)=0\}.
\]
Define
\[
 \mathcal M_H:\Lin_0(H(\A),\bfk)
 \longrightarrow G_{H(\A)}(\bfk),
 \qquad
 \beta\longmapsto
 \Ext\!\left(\exp^{\convH}(\beta)\circ\iota\right),
\]
where \(\exp^{\convH}\) and~\eqref{eq:tens-con} are analogous; it suffices to replace \(\deshat\), \(\convTbar\) and \(\varepsilon_{\overline{T}(\A)}\) by \(\deshah\), \(\convH\) and \(\varepsilon_{H(\A)}\) respectively.

\begin{theorem}\label{thm:expanded-diagram}
The following diagram is commutative.
\[
\begin{tikzcd}[column sep=huge,row sep=large]
\Lin_0(\overline T(\A),\bfk)
  \arrow[r,"c\mapsto c\circ\conc"]
  \arrow[d,"P^{\mathrm{cla}}"']
&
\Lin_0(H(\A),\bfk)
  \arrow[d,"\mathcal M_H"]
\\
\bfk^{\Part(I)}
  \arrow[d,"(\clNC)_*"']
&
G_{H(\A)}(\bfk)
  \arrow[d,"\log_{\prec}^{H(\A)}"]
\\
\bfk^{\NC(I)}
&
 g_{H(\A)}(\bfk)
  \arrow[l,"P^{\mathrm{fre}}"]
\end{tikzcd}
\]
\end{theorem}

\begin{proof}
Fix \(c\in\Lin_0(\overline T(\A),\bfk)\).
Since \(\conc\) is a Hopf algebra morphism by
Theorem~\ref{thm:flattening-morphism}~\Item ~\ref{item:conc-hopfsurhom}, precomposition
with \(\conc\) preserves convolution products and their units.
Consequently, for every integer \(q\geq0\),
\[
 (c\circ\conc)^{\convH q}
 =c^{\convTbar q}\circ\conc.
\]
Then 
\[
 \exp^{\convH}(c\circ\conc)
 =\bigl(\exp^{\convTbar}(c)\bigr)\circ\conc
 \overset{\eqref{eq:tens-con}}{=}\mu_c\circ\conc
 \overset{\eqref{eq:phic-con}}{=}\Phi_c^{\mathfrak c}.
\]
Since \(\conc\circ\iota=\idTbar\), it follows that
\[
 \mathcal M_H(c\circ\conc)
 =\Ext\!\left(\exp^{\convH}(\cconc)\circ\iota\right)
 =\Ext(\mu_c\circ\conc\circ\iota)
 =\Ext(\mu_c)
 \overset{\eqref{eq:muc-phic}}{=}\Phi_c.
\]
Thus the path through \(\Lin_0(H(\A),\bfk)\) and
\( G_{H(\A)}(\bfk)\) reaches
\(\log_{\prec}^{H(\A)}(\Phi_c)=r_c\in g_{H(\A)}(\bfk)\).
For every \(\pi\in\NC(I)\), 
\[
\begin{aligned}
&\Bigl(P^{\mathrm{fre}}
  \bigl(\log_{\prec}^{H(\A)}(\mathcal M_H(c\circ\conc))\bigr)\Bigr)(\pi)\\
 =&\ (r_c)_\pi\\
 =&\ \sum_{\substack{\sigma\in\Part(I)\\
                  \overline\sigma=\pi}}c_\sigma\hspace{1cm}(\text{by Lemma }\ref{lemma:closure-fiber-formula})\\
 =&\ \bigl((\clNC)_*\circ P^{\mathrm{cla}}(c)\bigr)(\pi).\hspace{1cm}(\text{by Definition }\ref{def:closure-pushforward})
\end{aligned}
\]
This holds for every \(c\) and \(\pi\), so the diagram commutes.
\end{proof}

\begin{remark}
\begin{enumerate}
\item Although Theorem~\ref{thm:expanded-diagram} is formulated for arbitrary
\(c\in\Lin_0(\overline T(\A),\bfk)\), it has a direct interpretation
in non-commutative probability.
For a non-commutative probability space $(\A,\varphi)$ with classical
cumulant functional $c$, $\Phi_c=\mathcal M_H(c\circ\conc)$ is its
moment character, and $r_c=\log_{\prec}^{H(\A)}(\Phi_c)$ is the
infinitesimal character encoding its free cumulants. Thus the
transform relates classical and free cumulants of the same
probability space.

\item
Our construction is \emph{algebraic and global}, acting on entire
cumulant functionals simultaneously in all word lengths.
Lehner's formula~\cite[Theorem~4.1 and its proof]{Lehner2002}
provides a \emph{combinatorial and local} description through
sums over individual non-crossing closure fibers.
Theorem~\ref{thm:expanded-diagram} identifies these functional-level
and fiberwise descriptions, connecting the closure-fiber formula
with our Hopf-algebraic construction and perturbed half-shuffle
equations.
\end{enumerate}
\end{remark}

\begin{example}\label{ex:degree-four-closure}
Take \(I=[4]\) and assume \(c_{\{i\}}=0\) for \(i\in I\).
We compute the two paths in Theorem~\ref{thm:expanded-diagram}
at the one-block partition \(\pi=\onepart{I}\).
\begin{enumerate}
\item \textbf{The functional path.}
\[
 c\longmapsto\beta:=c\circ\conc
 \xmapsto{\mathcal M_H}\Phi_c
 \xmapsto{\log_{\prec}^{H(\A)}}r_c
 \xmapsto{P^{\mathrm{fre}}}P^{\mathrm{fre}}(r_c).
\]
The first arrow gives \(\beta(a_S)=c_S\) and
\(\beta(a_I^S)=c_{I\setminus S}\); in particular,
\(\beta(a_2|a_4)=c_{\{2,4\}}\).
Since the singleton values vanish, the convolution square is
\[
\begin{aligned}
 \beta^{\convH 2}(a_I)
 =&\ \sum_{\substack{S\subset I\\|S|=2}}c_Sc_{I\setminus S}\\
 =&\ 2\bigl(c_{\{1,2\}}c_{\{3,4\}}
          +c_{\{1,3\}}c_{\{2,4\}}
          +c_{\{1,4\}}c_{\{2,3\}}\bigr).
\end{aligned}
\]
The powers \(\beta^{\convH q}(a_I)\) vanish for \(q\geq3\):
three or four nonempty parts contain a singleton, and higher powers
vanish by degree. Hence \(\Phi_c=\mathcal M_H(\beta)\) gives
\[
\begin{aligned}
 m_I
 =&\ \exp^{\convH}(\beta)(a_I)\\
 =&\ \beta(a_I)+\frac12\beta^{\convH 2}(a_I)\\
 =&\ c_I+c_{\{1,2\}}c_{\{3,4\}}
       +c_{\{1,3\}}c_{\{2,4\}}
       +c_{\{1,4\}}c_{\{2,3\}}.
\end{aligned}
\]
The same calculation gives \(m_{\{i\}}=0\) and
\(m_{\{i,j\}}=c_{\{i,j\}}\).
For the logarithm arrow,~\eqref{eq:left-half-logarithm} first yields
\(r_c(a_i)=0\) and \(r_c(a_ia_j)=c_{\{i,j\}}\) for \(i<j\).
For \(S=\{1,3\}\), bar-multiplicativity gives
\(\Phi_c(a_2|a_4)=m_{\{2\}}m_{\{4\}}=0\).
Thus, on \(a_I\), only \(S=I,\{1,2\},\{1,4\}\) can contribute
in that half-shuffle identity. The last two arrows give
\[
\begin{aligned}
 \bigl(P^{\mathrm{fre}}(r_c)\bigr)(\pi)
 =&\ r_c(a_I)\\
 =&\ m_I-c_{\{1,2\}}c_{\{3,4\}}
        -c_{\{1,4\}}c_{\{2,3\}}\\
 =&\ c_I+c_{\{1,3\}}c_{\{2,4\}}.
\end{aligned}
\]

\item \textbf{The path through \(\bfk^{\Part(I)}\).}
\[
 c\xmapsto{P^{\mathrm{cla}}}w:=P^{\mathrm{cla}}(c)
 \xmapsto{(\clNC)_*}(\clNC)_*w.
\]
Put \[\sigma_\times:=\{\{1,3\},\{2,4\}\}.\]
The first arrow gives
\[
 w(\pi)=c_I,
 \qquad w(\sigma_\times)=c_{\{1,3\}}c_{\{2,4\}}.
\]
Since \(\sigma_\times\) is the only crossing partition of \(I\)
and its closure is \(\pi\),
\[
 \clNC^{-1}(\pi)=\{\pi,\sigma_\times\}.
\]
The push-forward arrow therefore gives
\[
 \bigl((\clNC)_*w\bigr)(\pi)
 =w(\pi)+w(\sigma_\times)
 =c_I+c_{\{1,3\}}c_{\{2,4\}}.
\]
\end{enumerate}
For every other \(\rho\in\NC(I)\), both paths give \(c_\rho\):
the products \(c_\rho\) for the two non-crossing pairings are unchanged,
and \(c_\rho=0\) whenever \(\rho\) contains a singleton.
\end{example}

\appendix
\numberwithin{equation}{section}

\section{Low-degree calculations for the perturbed recursions}
\label{app:low-degree-calculations}

In this Appendix, we compute the two perturbed recursions in degrees
three, four, and five. The calculations illustrate how the moment-perturbed and cumulant-perturbed fixed-point equations give the classical
moment--cumulant expansions established in
Theorem~\ref{thm:classical-perturbations}.

We retain the ordered-word notation of
Subsection~\ref{subsec:construction-dcanc}, the moment and cumulant
coordinates of Subsection~\ref{subsec:left-moment-perturbation}, the complement-gap notation~\eqref{eq:aI-conv} and the multilinear moment coordinates~\eqref{eq:mom-cum-func}. 

\subsection{The moment-perturbed fixed-point equation}
\label{app:left-calculations}

We first apply the moment-perturbed fixed-point equation~\eqref{eq:classical-left-fixed-point}. In each calculation, we evaluate
the complement bar words by the perturbed moment functional and then
expand the remaining moments recursively.

\subsubsection{The third-order calculation}
Let
$
 A_3:=a_1a_2a_3.
$
The left half-coproduct has \(2^2=4\) terms, indexed by the
subsets of \([3]\) containing \(\min[3]=1\):
\[
 \DprecH(A_3)
 =a_1a_2a_3\otimes\mathbf1
 +a_1a_2\otimes a_3
 +a_1a_3\otimes a_2
 +a_1\otimes a_2a_3.
\]
The right-hand factors are evaluated by \(\Phiconc\) from
\eqref{eq:perturbed-moment-functional}, as shown below. In this
and the subsequent tables, digit strings abbreviate index sets, and commas between such
strings separate the connected components of a complement.
\begin{center}
\small
\begin{tabular}{c@{\qquad}c@{\qquad}c@{\qquad}c}
\toprule
\(S\)&\([3]\setminus S\)&\(a_{[3]}^S\)&\(\Phiconc(a_{[3]}^S)\)\\
\midrule
123&\(\varnothing\)&\(\mathbf1\)&1\\
12&3&\(a_3\)&\(m_{\{3\}}\)\\
13&2&\(a_2\)&\(m_{\{2\}}\)\\
1&23&\(a_2a_3\)&\(m_{\{2,3\}}\)\\
\bottomrule
\end{tabular}
\end{center}
Evaluating~\eqref{eq:classical-left-fixed-point} on the ordinary word
\(A_3\) gives
\[
 \Phi(A_3)
 =(c_H\prec\Phiconc)(A_3)
 =\bigl(c_H\prec(\Phi+\Phi\circ\ioconc)\bigr)(A_3),
\]
which gives
\[
 m_{\{1,2,3\}}
 =\Phi(a_1a_2a_3)
 =c_{\{1,2,3\}}
  +c_{\{1,2\}}m_{\{3\}}
  +c_{\{1,3\}}m_{\{2\}}
  +c_{\{1\}}m_{\{2,3\}}.
\]
The same recursion gives the lower-order moments
\[
 m_{\{2\}}=c_{\{2\}},
 \qquad
 m_{\{3\}}=c_{\{3\}},
 \qquad
 m_{\{2,3\}}=c_{\{2,3\}}+c_{\{2\}}c_{\{3\}}.
\]
Thus
\[
 m_{\{1,2,3\}}
 =c_{\{1,2,3\}}+
 c_{\{1\}}c_{\{2,3\}}
 +c_{\{1,2\}}c_{\{3\}}
 +c_{\{1,3\}}c_{\{2\}}
 +c_{\{1\}}c_{\{2\}}c_{\{3\}}.
\]
The \(1+3+1=5\) terms correspond to the five set partitions of \([3]\).

\subsubsection{The fourth-order calculation}
Let
$
 A_4:=a_1a_2a_3a_4.
$
The left half-coproduct has \(2^3=8\) terms, indexed by the
subsets of \([4]\) containing \(1\):
\[
 \begin{aligned}
 \DprecH(A_4)
 =&\ a_1a_2a_3a_4\otimes\mathbf1
 +a_1a_2a_3\otimes a_4
 +a_1a_2a_4\otimes a_3
 +a_1a_3a_4\otimes a_2\\
 &\ +a_1a_2\otimes a_3a_4
 +a_1a_3\otimes a_2|a_4
 +a_1a_4\otimes a_2a_3
 +a_1\otimes a_2a_3a_4.
 \end{aligned}
\]
Evaluating the right-hand factors by \(\Phiconc\) gives
\begin{center}
\small
\begin{tabular}{c@{\qquad}c@{\qquad}c@{\qquad}c}
\toprule
\(S\)&\([4]\setminus S\)&\(a_{[4]}^S\)&\(\Phiconc(a_{[4]}^S)\)\\
\midrule
1234&\(\varnothing\)&\(\mathbf1\)&1\\
123&4&\(a_4\)&\(m_{\{4\}}\)\\
124&3&\(a_3\)&\(m_{\{3\}}\)\\
134&2&\(a_2\)&\(m_{\{2\}}\)\\
12&34&\(a_3a_4\)&\(m_{\{3,4\}}\)\\
13&2,4&\(a_2|a_4\)&\(m_{\{2,4\}}\)\\
14&23&\(a_2a_3\)&\(m_{\{2,3\}}\)\\
1&234&\(a_2a_3a_4\)&\(m_{\{2,3,4\}}\)\\
\bottomrule
\end{tabular}
\end{center}
The only disconnected complement in this table corresponds to
\(S=\{1,3\}\):
\[
 [4]\setminus S=\{2,4\},
 \qquad a_{[4]}^S=a_2|a_4,
\]
and the perturbation changes its value from the unperturbed product
\(\Phi(a_2|a_4)=m_{\{2\}}m_{\{4\}}\) to the concatenated moment
\[
 \Phiconc(a_2|a_4)=\Phi(a_2a_4)=m_{\{2,4\}}.
\]
Substitution into the perturbed fixed-point equation gives
\[
 \begin{aligned}
 m_{\{1,2,3,4\}}
=&\ \Phi(a_1a_2a_3a_4)\\
=&\ c_{\{1,2,3,4\}}
 +c_{\{1,2,3\}}m_{\{4\}}
 +c_{\{1,2,4\}}m_{\{3\}}
 +c_{\{1,3,4\}}m_{\{2\}}\\
&\ +c_{\{1,2\}}m_{\{3,4\}}
 +c_{\{1,3\}}m_{\{2,4\}}
 +c_{\{1,4\}}m_{\{2,3\}}
 +c_{\{1\}}m_{\{2,3,4\}}.
 \end{aligned}
\]
The lower-order moments satisfy
\begin{align*}
 m_{\{3,4\}}=&\ c_{\{3,4\}}+c_{\{3\}}c_{\{4\}},\\
 m_{\{2,4\}}=&\ c_{\{2,4\}}+c_{\{2\}}c_{\{4\}},\\
 m_{\{2,3\}}=&\ c_{\{2,3\}}+c_{\{2\}}c_{\{3\}},
\end{align*}
and
\[
 m_{\{2,3,4\}}=c_{\{2,3,4\}}
  +c_{\{2\}}c_{\{3,4\}}
  +c_{\{2,3\}}c_{\{4\}}
  +c_{\{2,4\}}c_{\{3\}}
  +c_{\{2\}}c_{\{3\}}c_{\{4\}}.
\]
Substituting these identities gives
\[
\begin{aligned}
 m_{\{1,2,3,4\}}
=&\ c_{\{1,2,3,4\}}\\
&\ +\Big(
 c_{\{1\}}c_{\{2,3,4\}}
 +c_{\{1,2\}}c_{\{3,4\}}
 +c_{\{1,2,3\}}c_{\{4\}}
 +c_{\{1,2,4\}}c_{\{3\}}\\
&\ \qquad
 +c_{\{1,3\}}c_{\{2,4\}}
 +c_{\{1,3,4\}}c_{\{2\}}
 +c_{\{1,4\}}c_{\{2,3\}}
 \Big)\\
&\ +\Big(
 c_{\{1\}}c_{\{2\}}c_{\{3,4\}}
 +c_{\{1\}}c_{\{2,3\}}c_{\{4\}}
 +c_{\{1\}}c_{\{2,4\}}c_{\{3\}}\\
&\ \qquad
 +c_{\{1,2\}}c_{\{3\}}c_{\{4\}}
 +c_{\{1,3\}}c_{\{2\}}c_{\{4\}}
 +c_{\{1,4\}}c_{\{2\}}c_{\{3\}}
 \Big)\\
&\ +c_{\{1\}}c_{\{2\}}c_{\{3\}}c_{\{4\}}.
\end{aligned}
\]
The \(1+7+6+1=15\) terms correspond to the fifteen set partitions of
\([4]\). In particular, the crossing two-block partition
\(\{\{1,3\},\{2,4\}\}\) appears through the term \(c_{\{1,3\}}c_{\{2,4\}}\), which has no counterpart indexed by that partition in the free
moment--cumulant expansion.

\subsubsection{The fifth-order calculation}

Let \(A_5:=a_1a_2a_3a_4a_5\). By~\eqref{eq:half-coproducts-H}, the
left half-coproduct has \(2^4=16\) terms, indexed by subsets
\(S\subseteq[5]\) containing \(1\). Their right-hand factors are
evaluated by \(\Phiconc\), as follows.
\begin{center}
\small
\begin{tabular}{c@{\qquad}c@{\qquad}c@{\qquad}c}
\toprule
\(S\)&\([5]\setminus S\)&\(a_{[5]}^S\)&\(\Phiconc(a_{[5]}^S)\)\\
\midrule
12345&\(\varnothing\)&\(\mathbf1\)&1\\
1234&5&\(a_5\)&\(m_{\{5\}}\)\\
1235&4&\(a_4\)&\(m_{\{4\}}\)\\
1245&3&\(a_3\)&\(m_{\{3\}}\)\\
1345&2&\(a_2\)&\(m_{\{2\}}\)\\
123&45&\(a_4a_5\)&\(m_{\{4,5\}}\)\\
124&3,5&\(a_3|a_5\)&\(m_{\{3,5\}}\)\\
125&34&\(a_3a_4\)&\(m_{\{3,4\}}\)\\
134&2,5&\(a_2|a_5\)&\(m_{\{2,5\}}\)\\
135&2,4&\(a_2|a_4\)&\(m_{\{2,4\}}\)\\
145&23&\(a_2a_3\)&\(m_{\{2,3\}}\)\\
12&345&\(a_3a_4a_5\)&\(m_{\{3,4,5\}}\)\\
13&2,45&\(a_2|a_4a_5\)&\(m_{\{2,4,5\}}\)\\
14&23,5&\(a_2a_3|a_5\)&\(m_{\{2,3,5\}}\)\\
15&234&\(a_2a_3a_4\)&\(m_{\{2,3,4\}}\)\\
1&2345&\(a_2a_3a_4a_5\)&\(m_{\{2,3,4,5\}}\)\\
\bottomrule
\end{tabular}
\end{center}
The ordinary-word equation~\eqref{eq:classical-left-fixed-point}, with
its normalization at the unit, is
\begin{equation*}
\begin{aligned}
 \Phi\circ\iota
 =&\ \varepsilon_{\overline T(\A)}
    +(c_H\prec\Phiconc)\circ\iota\\
 =&\ \varepsilon_{\overline T(\A)}
    +\bigl(c_H\prec(\Phi+\Phi\circ\ioconc)\bigr)\circ\iota.
\end{aligned}
\end{equation*}
Evaluating it on \(A_5\) gives
\begin{align}
 m_{\{1,2,3,4,5\}}
=&\ \Phi(a_1a_2a_3a_4a_5)\notag\\
=&\ c_{\{1,2,3,4,5\}}
+c_{\{1,2,3,4\}}m_{\{5\}}
+c_{\{1,2,3,5\}}m_{\{4\}}
+c_{\{1,2,4,5\}}m_{\{3\}}
+c_{\{1,3,4,5\}}m_{\{2\}}\notag\\
&\ +c_{\{1,2,3\}}m_{\{4,5\}}
+c_{\{1,2,4\}}m_{\{3,5\}}
+c_{\{1,2,5\}}m_{\{3,4\}}\notag\\
&\ +c_{\{1,3,4\}}m_{\{2,5\}}
+c_{\{1,3,5\}}m_{\{2,4\}}
+c_{\{1,4,5\}}m_{\{2,3\}}\notag\\
&\ +c_{\{1,2\}}m_{\{3,4,5\}}
+c_{\{1,3\}}m_{\{2,4,5\}}
+c_{\{1,4\}}m_{\{2,3,5\}}\notag\\
&\ +c_{\{1,5\}}m_{\{2,3,4\}}
+c_{\{1\}}m_{\{2,3,4,5\}}.
\label{eq:classical-direct-five}
\end{align}
This is the fifth-order instance of the classical first-block recursion. Recursively expanding every
moment in~\eqref{eq:classical-direct-five} by
\eqref{eq:classical-first-block-recursion} gives all set partitions of
\([5]\). Hence there are
\[
 |\Part([5])|=52
\]
terms, grouped below by the number of blocks:
\begingroup
\small
\begin{equation*}
\begin{aligned}
&m_{\{1,2,3,4,5\}}\\
=&\ \Phi(a_1a_2a_3a_4a_5)\\
=&\ c_{\{1,2,3,4,5\}}\\
&\ +\Big(
 c_{\{1\}}c_{\{2,3,4,5\}}
+c_{\{1,2\}}c_{\{3,4,5\}}
+c_{\{1,2,3\}}c_{\{4,5\}}
+c_{\{1,2,3,4\}}c_{\{5\}}
+c_{\{1,2,3,5\}}c_{\{4\}}\\
&\ \qquad
+c_{\{1,2,4\}}c_{\{3,5\}}
+c_{\{1,2,4,5\}}c_{\{3\}}
+c_{\{1,2,5\}}c_{\{3,4\}}
+c_{\{1,3\}}c_{\{2,4,5\}}
+c_{\{1,3,4\}}c_{\{2,5\}}\\
&\ \qquad
+c_{\{1,3,4,5\}}c_{\{2\}}
+c_{\{1,3,5\}}c_{\{2,4\}}
+c_{\{1,4\}}c_{\{2,3,5\}}
+c_{\{1,4,5\}}c_{\{2,3\}}
+c_{\{1,5\}}c_{\{2,3,4\}}
\Big)\\
&\ +\Big(
 c_{\{1\}}c_{\{2\}}c_{\{3,4,5\}}
+c_{\{1\}}c_{\{2,3\}}c_{\{4,5\}}
+c_{\{1\}}c_{\{2,3,4\}}c_{\{5\}}
+c_{\{1\}}c_{\{2,3,5\}}c_{\{4\}}\\
&\ \qquad
+c_{\{1\}}c_{\{2,4\}}c_{\{3,5\}}
+c_{\{1\}}c_{\{2,4,5\}}c_{\{3\}}
+c_{\{1\}}c_{\{2,5\}}c_{\{3,4\}}
+c_{\{1,2\}}c_{\{3\}}c_{\{4,5\}}\\
&\ \qquad
+c_{\{1,2\}}c_{\{3,4\}}c_{\{5\}}
+c_{\{1,2\}}c_{\{3,5\}}c_{\{4\}}
+c_{\{1,2,3\}}c_{\{4\}}c_{\{5\}}
+c_{\{1,2,4\}}c_{\{3\}}c_{\{5\}}\\
&\ \qquad
+c_{\{1,2,5\}}c_{\{3\}}c_{\{4\}}
+c_{\{1,3\}}c_{\{2\}}c_{\{4,5\}}
+c_{\{1,3\}}c_{\{2,4\}}c_{\{5\}}
+c_{\{1,3\}}c_{\{2,5\}}c_{\{4\}}\\
&\ \qquad
+c_{\{1,3,4\}}c_{\{2\}}c_{\{5\}}
+c_{\{1,3,5\}}c_{\{2\}}c_{\{4\}}
+c_{\{1,4\}}c_{\{2\}}c_{\{3,5\}}
+c_{\{1,4\}}c_{\{2,3\}}c_{\{5\}}\\
&\ \qquad
+c_{\{1,4\}}c_{\{2,5\}}c_{\{3\}}
+c_{\{1,4,5\}}c_{\{2\}}c_{\{3\}}
+c_{\{1,5\}}c_{\{2\}}c_{\{3,4\}}
+c_{\{1,5\}}c_{\{2,3\}}c_{\{4\}}\\
&\ \qquad
+c_{\{1,5\}}c_{\{2,4\}}c_{\{3\}}
\Big)\\
&\ +\Big(
 c_{\{1\}}c_{\{2\}}c_{\{3\}}c_{\{4,5\}}
+c_{\{1\}}c_{\{2\}}c_{\{3,4\}}c_{\{5\}}
+c_{\{1\}}c_{\{2\}}c_{\{3,5\}}c_{\{4\}}\\
&\ \qquad
+c_{\{1\}}c_{\{2,3\}}c_{\{4\}}c_{\{5\}}
+c_{\{1\}}c_{\{2,4\}}c_{\{3\}}c_{\{5\}}
+c_{\{1\}}c_{\{2,5\}}c_{\{3\}}c_{\{4\}}
+c_{\{1,2\}}c_{\{3\}}c_{\{4\}}c_{\{5\}}\\
&\ \qquad
+c_{\{1,3\}}c_{\{2\}}c_{\{4\}}c_{\{5\}}
+c_{\{1,4\}}c_{\{2\}}c_{\{3\}}c_{\{5\}}
+c_{\{1,5\}}c_{\{2\}}c_{\{3\}}c_{\{4\}}
\Big)\\
&\ +c_{\{1\}}c_{\{2\}}c_{\{3\}}c_{\{4\}}c_{\{5\}}.
\end{aligned}
\end{equation*}
\endgroup
The groups contain \(1,15,25,10,1\) terms, respectively, and hence
account for all \(52\) set partitions.

\Needspace{7\baselineskip}
\subsection{The cumulant-perturbed fixed-point equation}
\label{app:right-calculations}

We next apply the cumulant-perturbed fixed-point equation~\eqref{eq:right-fixed-point}
to the same words. Here the complement bar words are evaluated by the
perturbed cumulant functional; the resulting moment expansions agree
with those in Appendix~\ref{app:left-calculations}.

\subsubsection{The third-order calculation}

For \(A_3=a_1a_2a_3\), the right half-coproduct has
\(2^2=4\) terms, indexed by the subsets of \(\{2,3\}\):
\[
 \DsuccH(A_3)
 =\mathbf1\otimes a_1a_2a_3
 +a_2\otimes a_1|a_3
 +a_3\otimes a_1a_2
 +a_2a_3\otimes a_1.
\]
The right factors are evaluated by \(\cconc\) from
\eqref{eq:perturbed-cumulant-functional}:
\begin{center}
\small
\begin{tabular}{c@{\qquad}c@{\qquad}c@{\qquad}c}
\toprule
\(S\)&\(a_S\)&\(a_{[3]}^S\)&\(\cconc(a_{[3]}^S)\)\\
\midrule
\(\varnothing\)&\(\mathbf1\)&\(a_1a_2a_3\)&\(c_{\{1,2,3\}}\)\\
2&\(a_2\)&\(a_1|a_3\)&\(c_{\{1,3\}}\)\\
3&\(a_3\)&\(a_1a_2\)&\(c_{\{1,2\}}\)\\
23&\(a_2a_3\)&\(a_1\)&\(c_{\{1\}}\)\\
\bottomrule
\end{tabular}
\end{center}
Substitution into~\eqref{eq:right-fixed-point} gives
\[
\begin{aligned}
 m_{\{1,2,3\}}
 =&\ \Phi(a_1a_2a_3)\\
 =&\ c_{\{1,2,3\}}
 +m_{\{2\}}c_{\{1,3\}}
 +m_{\{3\}}c_{\{1,2\}}
 +m_{\{2,3\}}c_{\{1\}}\\
 =&\ c_{\{1,2,3\}}
 +c_{\{1,3\}}m_{\{2\}}
 +c_{\{1,2\}}m_{\{3\}}
 +c_{\{1\}}m_{\{2,3\}}.
\end{aligned}
\]
The same recursion gives the lower-order moments
\[
 m_{\{2\}}=c_{\{2\}},
 \qquad
 m_{\{3\}}=c_{\{3\}},
 \qquad
 m_{\{2,3\}}=c_{\{2,3\}}+c_{\{2\}}c_{\{3\}}.
\]
Hence
\[
 m_{\{1,2,3\}}
=c_{\{1,2,3\}}
+c_{\{1\}}c_{\{2,3\}}
 +c_{\{1,2\}}c_{\{3\}}
 +c_{\{1,3\}}c_{\{2\}}
 +c_{\{1\}}c_{\{2\}}c_{\{3\}}.
\]
The \(1+3+1=5\) terms again correspond to the set partitions of \([3]\).

\subsubsection{The fourth-order calculation}

For \(A_4=a_1a_2a_3a_4\), the right half-coproduct has
\(2^3=8\) terms, indexed by the subsets of \(\{2,3,4\}\):
\[
\begin{aligned}
 \DsuccH(A_4)
=&\ \mathbf1\otimes a_1a_2a_3a_4
 +a_2\otimes a_1|a_3a_4
 +a_3\otimes a_1a_2|a_4
 +a_4\otimes a_1a_2a_3\\
&\ +a_2a_3\otimes a_1|a_4
 +a_2a_4\otimes a_1|a_3
 +a_3a_4\otimes a_1a_2
 +a_2a_3a_4\otimes a_1.
\end{aligned}
\]
Evaluating the right-hand factors by \(\cconc\) gives
\begin{center}
\small
\begin{tabular}{c@{\qquad}c@{\qquad}c@{\qquad}c}
\toprule
\(S\)&\(a_S\)&\(a_{[4]}^S\)&\(\cconc(a_{[4]}^S)\)\\
\midrule
\(\varnothing\)&\(\mathbf1\)&\(a_1a_2a_3a_4\)&\(c_{\{1,2,3,4\}}\)\\
2&\(a_2\)&\(a_1|a_3a_4\)&\(c_{\{1,3,4\}}\)\\
3&\(a_3\)&\(a_1a_2|a_4\)&\(c_{\{1,2,4\}}\)\\
4&\(a_4\)&\(a_1a_2a_3\)&\(c_{\{1,2,3\}}\)\\
23&\(a_2a_3\)&\(a_1|a_4\)&\(c_{\{1,4\}}\)\\
24&\(a_2a_4\)&\(a_1|a_3\)&\(c_{\{1,3\}}\)\\
34&\(a_3a_4\)&\(a_1a_2\)&\(c_{\{1,2\}}\)\\
234&\(a_2a_3a_4\)&\(a_1\)&\(c_{\{1\}}\)\\
\bottomrule
\end{tabular}
\end{center}
Substitution into~\eqref{eq:right-fixed-point} gives
\[
\begin{aligned}
 m_{\{1,2,3,4\}}
=&\ \Phi(a_1a_2a_3a_4)\\
=&\ c_{\{1,2,3,4\}}
 +m_{\{2\}}c_{\{1,3,4\}}
 +m_{\{3\}}c_{\{1,2,4\}}
 +m_{\{4\}}c_{\{1,2,3\}}\\
&\ +m_{\{2,3\}}c_{\{1,4\}}
 +m_{\{2,4\}}c_{\{1,3\}}
 +m_{\{3,4\}}c_{\{1,2\}}
 +m_{\{2,3,4\}}c_{\{1\}}\\
=&\ c_{\{1,2,3,4\}}
 +c_{\{1,2,3\}}m_{\{4\}}
 +c_{\{1,2,4\}}m_{\{3\}}
 +c_{\{1,3,4\}}m_{\{2\}}\\
&\ +c_{\{1,2\}}m_{\{3,4\}}
 +c_{\{1,3\}}m_{\{2,4\}}
 +c_{\{1,4\}}m_{\{2,3\}}
 +c_{\{1\}}m_{\{2,3,4\}}.
\end{aligned}
\]
Using the lower-order classical recursions for \(m_{\{2,3\}}\),
\(m_{\{2,4\}}\), \(m_{\{3,4\}}\), and \(m_{\{2,3,4\}}\), we obtain
\[
\begin{aligned}
 m_{\{1,2,3,4\}}
=&\ c_{\{1,2,3,4\}}\\
&\ +\Big(
 c_{\{1\}}c_{\{2,3,4\}}
 +c_{\{1,2\}}c_{\{3,4\}}
 +c_{\{1,2,3\}}c_{\{4\}}
 +c_{\{1,2,4\}}c_{\{3\}}\\
&\ \qquad
 +c_{\{1,3\}}c_{\{2,4\}}
 +c_{\{1,3,4\}}c_{\{2\}}
 +c_{\{1,4\}}c_{\{2,3\}}
 \Big)\\
&\ +\Big(
 c_{\{1\}}c_{\{2\}}c_{\{3,4\}}
 +c_{\{1\}}c_{\{2,3\}}c_{\{4\}}
 +c_{\{1\}}c_{\{2,4\}}c_{\{3\}}\\
&\ \qquad
 +c_{\{1,2\}}c_{\{3\}}c_{\{4\}}
 +c_{\{1,3\}}c_{\{2\}}c_{\{4\}}
 +c_{\{1,4\}}c_{\{2\}}c_{\{3\}}
 \Big)\\
&\ +c_{\{1\}}c_{\{2\}}c_{\{3\}}c_{\{4\}}.
\end{aligned}
\]
The groups contain \(1,7,6,1\) terms, respectively, giving all
\(15\) set partitions of \([4]\).

\subsubsection{The fifth-order calculation}

For \(A_5=a_1a_2a_3a_4a_5\), the right half-coproduct has
\(2^4=16\) terms, indexed by the subsets of \(\{2,3,4,5\}\).
The table records their factors and the corresponding cumulant values.
\begin{center}
\small
\begin{tabular}{c@{\qquad}c@{\qquad}c@{\qquad}c}
\toprule
\(S\)&\(a_S\)&\(a_{[5]}^S\)&\(\cconc(a_{[5]}^S)\)\\
\midrule
\(\varnothing\)&\(\mathbf1\)&\(a_1a_2a_3a_4a_5\)&\(c_{\{1,2,3,4,5\}}\)\\
2&\(a_2\)&\(a_1|a_3a_4a_5\)&\(c_{\{1,3,4,5\}}\)\\
3&\(a_3\)&\(a_1a_2|a_4a_5\)&\(c_{\{1,2,4,5\}}\)\\
4&\(a_4\)&\(a_1a_2a_3|a_5\)&\(c_{\{1,2,3,5\}}\)\\
5&\(a_5\)&\(a_1a_2a_3a_4\)&\(c_{\{1,2,3,4\}}\)\\
23&\(a_2a_3\)&\(a_1|a_4a_5\)&\(c_{\{1,4,5\}}\)\\
24&\(a_2a_4\)&\(a_1|a_3|a_5\)&\(c_{\{1,3,5\}}\)\\
25&\(a_2a_5\)&\(a_1|a_3a_4\)&\(c_{\{1,3,4\}}\)\\
34&\(a_3a_4\)&\(a_1a_2|a_5\)&\(c_{\{1,2,5\}}\)\\
35&\(a_3a_5\)&\(a_1a_2|a_4\)&\(c_{\{1,2,4\}}\)\\
45&\(a_4a_5\)&\(a_1a_2a_3\)&\(c_{\{1,2,3\}}\)\\
234&\(a_2a_3a_4\)&\(a_1|a_5\)&\(c_{\{1,5\}}\)\\
235&\(a_2a_3a_5\)&\(a_1|a_4\)&\(c_{\{1,4\}}\)\\
245&\(a_2a_4a_5\)&\(a_1|a_3\)&\(c_{\{1,3\}}\)\\
345&\(a_3a_4a_5\)&\(a_1a_2\)&\(c_{\{1,2\}}\)\\
2345&\(a_2a_3a_4a_5\)&\(a_1\)&\(c_{\{1\}}\)\\
\bottomrule
\end{tabular}
\end{center}
Hence, 
\begin{align}
 m_{\{1,2,3,4,5\}}
=&\ \Phi(a_1a_2a_3a_4a_5)\notag\\
=&\ c_{\{1,2,3,4,5\}}
+c_{\{1,2,3,4\}}m_{\{5\}}
+c_{\{1,2,3,5\}}m_{\{4\}}
+c_{\{1,2,4,5\}}m_{\{3\}}
+c_{\{1,3,4,5\}}m_{\{2\}}\notag\\
&\ +c_{\{1,2,3\}}m_{\{4,5\}}
+c_{\{1,2,4\}}m_{\{3,5\}}
+c_{\{1,2,5\}}m_{\{3,4\}}\notag\\
&\ +c_{\{1,3,4\}}m_{\{2,5\}}
+c_{\{1,3,5\}}m_{\{2,4\}}
+c_{\{1,4,5\}}m_{\{2,3\}}\notag\\
&\ +c_{\{1,2\}}m_{\{3,4,5\}}
+c_{\{1,3\}}m_{\{2,4,5\}}
+c_{\{1,4\}}m_{\{2,3,5\}}\notag\\
&\ +c_{\{1,5\}}m_{\{2,3,4\}}
+c_{\{1\}}m_{\{2,3,4,5\}}.
\label{eq:right-direct-five}
\end{align}
Thus~\eqref{eq:right-direct-five} agrees with
\eqref{eq:classical-direct-five}. The cumulant-perturbed fixed-point equation
gives the same fifth-order classical first-block recursion as the moment-perturbed fixed-point equation.

Recursively expanding all moments on the right-hand side of
\eqref{eq:right-direct-five} by~\eqref{eq:classical-first-block-recursion} gives
all set partitions of \([5]\). Grouping the terms by their number of
blocks gives
\begingroup
\small
\begin{align*}
m_{\{1,2,3,4,5\}}
=&\ \Phi(a_1a_2a_3a_4a_5)\\*
=&\ c_{\{1,2,3,4,5\}}\\*
&\ +\Big(
 c_{\{1\}}c_{\{2,3,4,5\}}
+c_{\{1,2\}}c_{\{3,4,5\}}
+c_{\{1,2,3\}}c_{\{4,5\}}
+c_{\{1,2,3,4\}}c_{\{5\}}
+c_{\{1,2,3,5\}}c_{\{4\}}\\*
&\ \qquad
+c_{\{1,2,4\}}c_{\{3,5\}}
+c_{\{1,2,4,5\}}c_{\{3\}}
+c_{\{1,2,5\}}c_{\{3,4\}}
+c_{\{1,3\}}c_{\{2,4,5\}}
+c_{\{1,3,4\}}c_{\{2,5\}}\\*
&\ \qquad
+c_{\{1,3,4,5\}}c_{\{2\}}
+c_{\{1,3,5\}}c_{\{2,4\}}
+c_{\{1,4\}}c_{\{2,3,5\}}
+c_{\{1,4,5\}}c_{\{2,3\}}
+c_{\{1,5\}}c_{\{2,3,4\}}
\Big)\\
&\ +\Big(
 c_{\{1\}}c_{\{2\}}c_{\{3,4,5\}}
+c_{\{1\}}c_{\{2,3\}}c_{\{4,5\}}
+c_{\{1\}}c_{\{2,3,4\}}c_{\{5\}}
+c_{\{1\}}c_{\{2,3,5\}}c_{\{4\}}\\*
&\ \qquad
+c_{\{1\}}c_{\{2,4\}}c_{\{3,5\}}
+c_{\{1\}}c_{\{2,4,5\}}c_{\{3\}}
+c_{\{1\}}c_{\{2,5\}}c_{\{3,4\}}
+c_{\{1,2\}}c_{\{3\}}c_{\{4,5\}}\\*
&\ \qquad
+c_{\{1,2\}}c_{\{3,4\}}c_{\{5\}}
+c_{\{1,2\}}c_{\{3,5\}}c_{\{4\}}
+c_{\{1,2,3\}}c_{\{4\}}c_{\{5\}}
+c_{\{1,2,4\}}c_{\{3\}}c_{\{5\}}\\*
&\ \qquad
+c_{\{1,2,5\}}c_{\{3\}}c_{\{4\}}
+c_{\{1,3\}}c_{\{2\}}c_{\{4,5\}}
+c_{\{1,3\}}c_{\{2,4\}}c_{\{5\}}
+c_{\{1,3\}}c_{\{2,5\}}c_{\{4\}}\\*
&\ \qquad
+c_{\{1,3,4\}}c_{\{2\}}c_{\{5\}}
+c_{\{1,3,5\}}c_{\{2\}}c_{\{4\}}
+c_{\{1,4\}}c_{\{2\}}c_{\{3,5\}}
+c_{\{1,4\}}c_{\{2,3\}}c_{\{5\}}\\*
&\ \qquad
+c_{\{1,4\}}c_{\{2,5\}}c_{\{3\}}
+c_{\{1,4,5\}}c_{\{2\}}c_{\{3\}}
+c_{\{1,5\}}c_{\{2\}}c_{\{3,4\}}
+c_{\{1,5\}}c_{\{2,3\}}c_{\{4\}}\\*
&\ \qquad
+c_{\{1,5\}}c_{\{2,4\}}c_{\{3\}}
\Big)\\
&\ +\Big(
 c_{\{1\}}c_{\{2\}}c_{\{3\}}c_{\{4,5\}}
+c_{\{1\}}c_{\{2\}}c_{\{3,4\}}c_{\{5\}}
+c_{\{1\}}c_{\{2\}}c_{\{3,5\}}c_{\{4\}}\\*
&\ \qquad
+c_{\{1\}}c_{\{2,3\}}c_{\{4\}}c_{\{5\}}
+c_{\{1\}}c_{\{2,4\}}c_{\{3\}}c_{\{5\}}
+c_{\{1\}}c_{\{2,5\}}c_{\{3\}}c_{\{4\}}
+c_{\{1,2\}}c_{\{3\}}c_{\{4\}}c_{\{5\}}\\*
&\ \qquad
+c_{\{1,3\}}c_{\{2\}}c_{\{4\}}c_{\{5\}}
+c_{\{1,4\}}c_{\{2\}}c_{\{3\}}c_{\{5\}}
+c_{\{1,5\}}c_{\{2\}}c_{\{3\}}c_{\{4\}}
\Big)\\
&\ +c_{\{1\}}c_{\{2\}}c_{\{3\}}c_{\{4\}}c_{\{5\}}.
\end{align*}
\endgroup
The groups contain \(1,15,25,10,1\) terms, respectively. Their sum is
\(52\), the number of set partitions of \([5]\), in agreement with
Theorem~\ref{thm:classical-perturbations}.

\vskip 0.2in

\noindent
{\bf Acknowledgments.} This work is supported by the National Natural Science Foundation of China (12571019), the Natural Science Foundation of Gansu Province (25JRRA644) and Innovative Fundamental Research Group Project of Gansu Province (23JRRA684).  

\noindent
{\bf Declaration of interests. } The authors have no conflicts of interest to disclose.

\noindent
{\bf Data availability. } Data sharing is not applicable as no new data were created or analyzed.

\end{document}